\documentclass{article}
\author{R. Wm. Gosper 
\thanks{rwg@osots.com}
\and Rich Schroeppel
\thanks{rschroe@sandia.gov
Sandia is a multiprogram laboratory operated by Sandia Corporation,
a Lockheed Martin Company, for the United States Department of Energy
under Contract DE-AC04-94AL85000.
}
}
\title{Somos Sequence Near-Addition Formulas \\ and Modular Theta Functions}

\begin{document}
\bibliographystyle{plain}

\maketitle

\abstract {
We have discovered conjectural near-addition formulas for Somos sequences.
We have preliminary evidence suggesting the existence of modular theta functions.
}

\section{Introduction}

In 2001, as part of a research project investigating more
efficient public
key cryptography (PKC), Rich Schroeppel asked Bill Gosper to look for
Somos sequence addition formulas.
Gosper found some very interesting results
 immediately, and
further developments continued through 2003.
Cheryl Beaver and Schroeppel
 also investigated modular versions of the dilogarithm function ${\rm Li}_2(x)$
and the Trilogarithm ${\rm Li}_3(x)$ \cite{RCS},
and Schroeppel did some preliminary work on modular theta functions.

Somos sequences and theta functions are both promising
approaches for use in public key cryptography.
Cryptographic applications of Somos sequences are explored in \cite{OSG}.
Our dilogarithm and trilogarithm results are interesting,
but it's not obvious how to apply them to PKC problems.

Further development of modular versions of the
special functions of numerical analysis seems possible.
Likely candidates are the error function, the logarithmic integral,
the gamma function, and
 perhaps Bessel functions
and hypergeometric functions.

The main stumbling block is that inequalities, limiting processes, and infinite
series are unavailable, and we must fall back on functional equations for
much of the work.  Formal differentiation sometimes works.
Functional equations are very limited for the error function, but
more variety is available for the other special functions.
We have not explored the p-adic possibilities, which might permit the
 reintroduction of some of the forbidden concepts.

Both Somos sequences and theta functions have near-addition formulas:
equations that relate $f(x+y) f(x-y)$
 to $f(x)$ and $f(y)$ and $f$ of nearby $x$ and $y$ values.
These can be used with the well-known double-and-add method to calculate
function values at large multiples of the argument.

Gosper's results on near-addition formulas for Somos sequences are
reported in section 2.
Schroeppel has been able to prove a few of the formulas,
including two of the determinant identities for Somos4 \cite{OSG}.
Section 3 details our brief excursion into modular theta functions.
There is no conclusion: this research seems open ended.

\subsection{Somos Sequence Background}

RS first learned of Somos sequences from Michael Somos, around 1988.
The Somos sequence of order $N$ begins with a block of $N$ 1s,
and is generated from a simple
non-linear recurrence.
For Somos4, the recurrence is $ a_k a_{k+4} = a_{k+1} a_{k+3} + a_{k+2}^2 $;
for Somos5, it is  $a_k a_{k+5} = a_{k+1} a_{k+4} + a_{k+2} a_{k+3}$;
Somos6 and 7 follow the patterns for 4 and 5, with more terms in the folded
dot product.
The recurrence can also be used to extend the sequences in the negative direction; they
are palindromic.
Somos4 begins 1,1,1,1,2,... .
One surprising property is that all the terms in Somos4 - Somos7
are integers.  This was discovered
 back in the 1940s
by Morgan Ward \cite{Ward1948, WardAMS1948}
 for a sequence including Somos4.
He called his sequences Elliptic Divisibility Sequences; Somos4
is the odd numbered terms from a particular EDS.
RS was at an MSRI number theory workshop shortly after learning of the sequences,
and the group spent some time trying to prove integrality.
Eventually Dean Hickerson and Janice Malouf independently proved that Somos6
is integral.  Our experimenting showed that
we could modify the sequence initial values in various ways while
apparently keeping the integrality property.
Also, introducing algebraic
coefficients into the equation, such as
$ a_k a_{k+4} = x a_{k+1} a_{k+3} + y a_{k+2}^2$, often produced polynomials
with integer coefficients, rather than the expected ensemble of rational functions.
(Of course, this doesn't matter much when the values are interpreted mod $P$,
or in a finite field, which can handle fractions just fine.)
The raw integer values of the sequences seem to grow roughly as $~C^{K^2}$.
The Somos4 and Somos5 sequences have a close connection with elliptic curves
and classical theta functions.
 The higher order Somos sequences may be connected to hyperelliptic curves.
There's a moderate amount of background material scattered
around the net; Jim Propp's Somos page \cite{Sloane}, and Sloane's sequence
database \cite{Sloane2}, and Zagier's problem 5 \cite{Zagier}, have useful material.
Background on theta functions is available in the Abramowitz \& Stegun
\cite{AbSteg} compendium of
special functions, now available on the net (as a scanned photocopy) at \cite{AbStegWeb}.


\section {Somos Sequence Near-Addition Formulae}

\def\mod{\mathop{\rm mod} \nolimits}
\def\twoline#1#2#3{\vbox{\hbox to\displaywidth{$\displaystyle{#1}$\hfill}
\vskip#2\hbox to\displaywidth{\hfill$\displaystyle{#3}$\hfill}}}
\def\Tilde{\char126\relax}
\font\eightrm=cmr8 \font\sixrm=cmr6 
\font\eighttt=cmtt8
\parindent=0pt\overfullrule=0pt
\parskip=5pt

{\bf Summary:}  Somos($nx$) is calculable in $O(\log n)$ time from three values
near Somos($x$), at least for orders 4 and 5.  Orders 6 and 7 require longer intervals
of values.  Along the way, we find addition formul\ae\ for Somos and Somos-like sequences
of polynomials and algebraics, and reduce some fifth order recurrences to fourth and
third order.  We find three-term, four-variable relations for most of these, as well
as for ordinary $\vartheta$ functions.  A sequence of polynomials obeying the Somos4
recurrence has a particularly
nice doubling formula.  Many of these results fall out of a very general
determinant identity.  For certain algebraic ``Somos'' sequences, we find closed
forms in terms of Chebychev polynomials.

{\bf Definitions:}  $a_n:=$ Somos4, $b_n:=$ Somos5, $\ldots, e_n:=$ Somos8, {\it i.e.},

$$\vbox{\baselineskip=30pt
\halign{$\hfill\quad#_n:=\;$&$\displaystyle#\hfil$&$\,\hfil=\;\;#\;\;$&$=\;\;\ldots,#,\ldots\hfill$\cr
a&{a_{n-1}a_{n-3}+a_{n-2}^2\over a_{n-4}}&a_{\,3-n}&2,1,1,1,1,2,3,7,23,59,314,1529\cr
b&{b_{n-1}b_{n-4}+b_{n-2}b_{n-3}\over b_{n-5}}&b_{\,4-n}&2,1,1,1,1,1,2,3,5,11,37,83,274\cr
c&{c_{n-1}c_{n-5}+c_{n-2}c_{n-4}+c_{n-3}^2\over c_{n-6}}&c_{\,5-n}&3,1,1,1,1,1,1,3,5,9,23,75,421\cr
d&{d_{n-1}d_{n-6}+d_{n-2}d_{n-5}+d_{n-3}d_{n-4}\over d_{\,n-7}}&d_{6-n}&3,1,1,1,1,1,1,1,3,5,9,17,41,137\cr
e&\cdots\cr}
}$$

(all appropriately palindromic) where the tabulated values start with subscript $n=-1$
to show the center of symmetry.
They are integer sequences until

{\bf Somos8:}  $\displaystyle e_{17}={420514\over7}$, so we're not too interested in Somos8.
On the other hand, Somos6 satisfies
$$ c_{n}={{-c_{n-1}c_{n-8}-c_{n-2}c_{n-7}+c_{n-3}c_{n-6}+34c_{n-4}c%
 _{n-5}}\over{c_{n-9}}}, $$
which is pretty much a Somos9.  And, for all $t$ and $u$, Somos4 satisfies
$$ a_{n}={{\left(t-7\right) a_{n-1} a_{n-7}+\left(u-5 t+31\right) a_{n-2} a_{%
 n-6}+\left(4 t-u+1\right) a_{n-3} a_{n-5}-u a^{2}_{n-4}}\over{a_{n-8}}}, $$
a double continuum of quasi-Somos8s.

Furthermore, for all $t$, the sequence $s_n:=a^2_n$ satisfies
$$ s_{n}={{\left(6-t\right) s_{n-1} s_{n-7}+\left(5 t-130\right) s_{n-2} s_{n%
 -6}+\left(749-4 t\right) s_{n-3} s_{n-5}+\left(20 t-4\right)%
  s^{2}_{n-4}}\over{s_{n-8}}}. $$

{\bf Change of variable:}
A Somos sequence may be multiplied by any constant.
A Somos sequence multiplied by an arbitrary geometric progression satisfies the
same recurrence, but usually loses its palindrome symmetry.
The ``odd'' (Somos5 and Somos7) sequences may also be termwise multiplied by any
number of factors of the form $\tan(x+n\pi/2)$ without even disturbing the
palindrome property.

There is, however, no sharp dichotomy between odd and even, since Somos4 satisfies
the quasiSomos5 (Quasimodo?) (odd) recurrence
$$ a_{n}={{5 a_{n-3} a_{n-2}-a_{n-4} a_{n-1}}\over{a_{n-5}}}, $$

as does $a_n\tan(x+\pi n/2)$, etc.

The sequence $s_n:= r^{n^2}a_n$ satisfies
$$ s_{n}={{ r^{6} s_{n-1}s_{n-3}+ r^{8}s^{2}_{n-2}}\over{s_{n-4}}},$$
while
$s_n:=r^{n^2}b_n$ satisfies
$$ s_{n}={{r^{8} s_{n-1} s_{n-4}+r^{12} s_{n-2} s_{n-3}}\over{s_{n-5}}}. $$

Similarly, $s_n:=r^{n^2}c_n$ satisfies
$$ s_{n}={{r^{10} s_{n-1} s_{n-5}+r^{16} s_{n-2} s_{n-4}+r^{18} s^{2}%
 _{n-3}}\over{s_{n-6}}}, $$
while $s_n:=r^{n^2}d_n$ satisfies
$$ s_{n}={{r^{12} s_{n-1} s_{n-6}+r^{20} s_{n-2} s_{n-5}+r^{24} s_{n-3} s_{n-%
 4}}\over{s_{n-7}}}. $$

{\bf Notation:}

\begin{eqnarray*}
D_s\pmatrix{x_1,&x_2,&\ldots,&x_n\cr y_1,&y_2,&\ldots,&y_n\cr}:&= &
    \det\left[s_{x_{i}-y_{j}}s_{x_{i}+y_{j}}\right]_{1\le i,j\le n}\cr
  &=&\left|\matrix{s_{x_1-y_1}s_{x_1+y_1}&s_{x_1-y_2}s_{x_1+y_2}&\ldots&
                    s_{x_1-y_n}s_{x_1+y_n}\cr
                  s_{x_2-y_1}s_{x_2+y_1}&s_{x_2-y_2}s_{x_2+y_2}&\ldots&
                    s_{x_2-y_n}s_{x_2+y_n}\cr
                  \vdots&\vdots&&\vdots\cr
                  s_{x_n-y_1}s_{x_n+y_1}&s_{x_n-y_2}s_{x_n+y_2}&\ldots&
                    s_{x_n-y_n}s_{x_n+y_n}\cr}\right|.
\end{eqnarray*}

Note that each term of the expanded determinant will have subscripts summing to
$2x_1+2x_2+\ldots+2x_n$.  This is decidedly not symmetrical in $x$ and $y$,
so that an identity involving a $D$ operator may yield a new identity under
interchange of the $x$ and $y$ vectors.

{\bf Conjecture 4:}  the determinant

$$D_a\pmatrix{u,&v,&w\cr x,&y,&z\cr}=
 \left|\matrix{a_{u-x}\,a_{u+x}&a_{u-y}\,a_{u+y}&a_{u-z}\,a_{u+z}%
 \cr a_{v-x}\,a_{v+x}&a_{v-y}\,a_{v+y}&a_{v-z}\,a_{v+z}\cr a_{%
 w-x}\,a_{w+x}&a_{w-y}\,a_{w+y}&a_{w-z}\,a_{w+z}\cr }\right|=0, $$

where $u,v,w,x,y,$ and $z$ are arbitrary integers.  {\it E.g.},

$$D_a \pmatrix{n-2,&0,&1\cr 0,&1,&2\cr } 
=\left| \matrix{a^{2}_{n-2}&a_{n-3}\,a_{n-1}&a_{n-4}\,a_{n}\cr 1&2&3\cr 1&1&2\cr }\right|
=-a_{n-4}\,a_{n}+a_{n-3}\,a_{n-1}+a^{2}_{n-2},$$

the defining recurrence for Somos4.

Note that the determinant also
vanishes for $a_t:=\sin t$, for arbitrary complex $u,v,w,x,y,$ and $z$.  More
interestingly, experimental Taylor expansion at $q=0$ plus several numerical
experiments suggest that the same
goes for $a_t:=\vartheta_j(t,q)$.  The published addition formul\ae\ mixing two or
more different $j$ are merely the result of choosing $v,w,y,$ and $z$ to be things
like $\pi/2$ and $\pi\tau/2$ (and $0$).  (Whittaker \& Watson, crediting  Jacobi,
 list numerous special cases, suggesting that the more general formula was not
 yet known.)

Still more generally,
$$ 0=\left|\matrix{\vartheta_{s}\left(x-u,q\right)\,\vartheta_{t}\left(x+u,q%
 \right)&\vartheta_{s}\left(y-u,q\right)\,\vartheta_{t}\left(y+u,q%
 \right)&\vartheta_{s}\left(z-u,q\right)\,\vartheta_{t}\left(z+u,q%
 \right)\cr \vartheta_{s}\left(x-v,q\right)\,\vartheta_{t}\left(x+v,q%
 \right)&\vartheta_{s}\left(y-v,q\right)\,\vartheta_{t}\left(y+v,q%
 \right)&\vartheta_{s}\left(z-v,q\right)\,\vartheta_{t}\left(z+v,q%
 \right)\cr \vartheta_{s}\left(x-w,q\right)\,\vartheta_{t}\left(x+w,q%
 \right)&\vartheta_{s}\left(y-w,q\right)\,\vartheta_{t}\left(y+w,q%
 \right)&\vartheta_{s}\left(z-w,q\right)\,\vartheta_{t}\left(z+w,q%
 \right)\cr }\right|. $$

{\it E.g.}, putting $s=1, u=x, v=y$,
\begin{eqnarray*}0\quad=&\;\;\;\;\vartheta_{1}\left(x-w,q\right)\,\vartheta_{t}\left(x+w,q\right)\,%
 \vartheta_{1}\left(z-y,q\right)\,\vartheta_{t}\left(z+y,q\right)\;\cr
 &-\;\vartheta_{1}\left(y-w,q\right)\,\vartheta_{t}\left(y+w,q\right)\,%
 \vartheta_{1}\left(z-x,q\right)\,\vartheta_{t}\left(z+x,q\right)\; &\hspace{1in}{\rm (4vars)}\cr
 &+\;\vartheta_{1}\left(y-x,q\right)\,\vartheta_{t}\left(y+x,q\right)\,%
 \vartheta_{1}\left(z-w,q\right)\,\vartheta_{t}\left(z+w,q\right),&\cr
\end{eqnarray*}
a three term identity in four variables.

{\bf Conjecture 4.5:}  The determinant

$$D_s\pmatrix{u+1/2,&v+1/2,&w+1/2\cr x-1/2,&y-1/2,&z-1/2\cr}=
 \left|\matrix{s_{1+u-x}\,s_{u+x}&s_{1+u-y}\,s_{u+y}&s_{1+u-z}\,s_{u+z}%
 \cr s_{1+v-x}\,s_{v+x}&s_{1+v-y}\,s_{v+y}&s_{1+v-z}\,s_{v+z}\cr s_{1%
 +w-x}\,s_{w+x}&s_{1+w-y}\,s_{w+y}&s_{1+w-z}\,s_{w+z}\cr }\right|=0, $$

where $s_n:=a_n$ or $b_n$, and $u,v,w,x,y,$ and $z$ are arbitrary integers.  {\it E.g.},
\begin{eqnarray*} 0&=&D_b\pmatrix{n-\mathchoice {{5}\over{2}}{{5}\over{2}}{5/2}{5/2},&%
 \mathchoice {{1}\over{2}}{{1}\over{2}}{1/2}{1/2},&\mathchoice {{3%
 }\over{2}}{{3}\over{2}}{3/2}{3/2}\cr
\noalign{\vskip3pt}
 \mathchoice {{1}\over{2}}{{1%
 }\over{2}}{1/2}{1/2},&\mathchoice {{3}\over{2}}{{3}\over{2}}{3/2}{3/2%
 },&\mathchoice {{5}\over{2}}{{5}\over{2}}{5/2}{5/2}\cr } \cr
\noalign{\vskip9pt}
&=&\left|\matrix{b_{n-3}\,b_{n-2}&b_{n-4}\,b_{n-1}&b_{n-5}\,b_{n}\cr 1&2&3%
 \cr 1&1&2\cr }\right|\;\;=\;\;-b_{n-5}\,b_{n}+b_{n-4}\,b_{n-1}+b_{n-3}\,b_{n-2},\cr
\end{eqnarray*}
the defining recurrence for Somos5.

{\bf Somos4 addition formul\ae:}
(See the section ``Somos4oid polynomials'' for a sequence $s_n$ with much nicer
 addition formul\ae\ than those derived here for $a_n$.)

Suppose we have four consecutive values $a_{x-1}, a_x, a_{x+1}, a_{x+2}$.  Choose $s:=a, u=x, y=0, z=-1, v=0, w=1$ to get
\begin{eqnarray*}D_a\pmatrix{x+1,&0,&1\cr x-1,&0,&1}&=&D_a\pmatrix{x+1/2,&1/2,&3/2\cr x-1/2,&1/2,&-3/2}\\
& & \\
&=& \left|\matrix{a_{1}\,a_{2\,x}&a_{x}\,a_{x+1}&a_{x-1}\,a_{x+2}\cr a_{1-x%
 }\,a_{x}&a_{0}\,a_{1}&a_{-1}\,a_{2}\cr a_{2-x}\,a_{x+1}&a_{1}\,a_{2}%
 &a_{0}\,a_{3}\cr }\right|\\
\noalign{\vskip6pt}
 &=&\left|\matrix{a_{2\,x}&a_{x}\,a_{x+1}&a_{x-1}\,a_{x+2}%
 \cr a_{x}\,a_{x+2}&1&2\cr a^{2}_{x+1}&1&1\cr }\right|\;=\;0\rm,\\
\end{eqnarray*}
giving us $a_{2x}$.  Alternatively,
\begin{eqnarray*}D_a\pmatrix{x+1/2,&1/2,&3/2\cr x-3/2,&3/2,&-1/2}&=&
  \left|\matrix{a_{2}\,a_{2\,x-1}&a_{x-1}\,a_{x+2}&a_{x}\,a_{x+1}\cr a_{2%
 -x}\,a_{x-1}&a_{-1}\,a_{2}&a_{0}\,a_{1}\cr a_{3-x}\,a_{x}&a_{0}\,a_{%
 3}&a_{1}\,a_{2}\cr }\right|\\
\noalign{\vskip5pt}
 &=&\left|\matrix{a_{2\,x-1}&a_{x-1}\,a_{x+2}&a_{x}\,a_{%
 x+1}\cr a_{x-1}\,a_{x+1}&2&1\cr a^{2}_{x}&1&1\cr }\right|\;=\;0,\\
\end{eqnarray*}
giving us $a_{2x-1}$.

Now run the Somos4 recurrence one step forward to get
$a_{x+3}$ and replace $x$ by $x+1$ in the preceding two
determinants to get the four consecutive values $a_{2x-1},a_{2x},a_{2x+1},a_{2x+2}$.  So we can double $a_{nx}$ to $a_{2nx}$.

Now suppose that we have the four values around $a_{nx}$ and also around $a_{(n+1)x}$.  Then 
\begin{eqnarray*}D_a\pmatrix{nx+1,&0,&1\cr(n+1)x-1,&0,&1}&=&D_a\pmatrix{(n+1)x+1/2,&1/2,&3/2\cr nx-1/2,&-3/2,&-1/2}\cr
\noalign{\vskip3pt}
&=&\left|\matrix{a_{x+1}\,a_{\left(2\,n+1\right)\,x}&a_{\left(n+1\right)\,%
 x-1}\,a_{\left(n+1\right)\,x+2}&a_{\left(n+1\right)\,x}\,a_{\left(n+%
 1\right)\,x+1}\cr a_{n\,x}\,a_{1-n\,x}&a_{-1}\,a_{2}&a_{0}\,a_{1}%
 \cr a_{2-n\,x}\,a_{n\,x+1}&a_{0}\,a_{3}&a_{1}\,a_{2}\cr }\right|\cr
\noalign{\vskip6pt}
&=&\left|\matrix{a%
 _{x+1}\,a_{\left(2\,n+1\right)\,x}&a_{\left(n+1\right)\,x-1}\,a_{%
 \left(n+1\right)\,x+2}&a_{\left(n+1\right)\,x}\,a_{\left(n+1\right)%
 \,x+1}\cr a_{n\,x}\,a_{n\,x+2}&2&1\cr a^{2}_{n\,x+1}&1&1\cr }\right|\;=\;0\rm.\cr
\end{eqnarray*}

So from $a_{nx}$ and $a_{(n+1)x}$ we get $a_{2nx}$ and $a_{(2n+1)x}$.  Thus we can multiply by maintaining eight values.  {\it E.g.},
$$105x \leftarrow
(53x, 52x)\leftarrow(27x,26x)\leftarrow(14x,13x)\leftarrow(7x,6x)\leftarrow(4x,3x)\leftarrow(2x,x)\rm.$$

In principle, we need only maintain two sets of three values,
$a_{nx-1},  a_{nx},  a_{nx+1},$ and \\
$a_{(n+1)x-1},  a_{(n+1)x},  a_{(n+1)x+1}$,
by virtue of the third order relation

{\bf Conjecture 4a} (``derived'' below){\bf:}
$$ a^{2}_{x-1}\,a^{2}_{x+2}+a^{3}_{x}\,a_{x+2}+a_{x-1}\,a^{3}_{x+1}+a%
 ^{2}_{x}\,a^{2}_{x+1}=4\,a_{x-1}\,a_{x}\,a_{x+1}\,a_{x+2}, $$
with which we can eliminate $a_{x+2}$ from:

$$ \left|\matrix{a_{2\,x-1}&a_{x-1}\,a_{x+2}&a_{x}\,a_{x+1}\cr a_{x-1}\,a%
 _{x+1}&2&1\cr a^{2}_{x}&1&1\cr }\right|=a_{2\,x-1}-a_{x-1}\,\left(a_{x-1}\,%
 a_{x+1}-a^{2}_{x}\right)\,a_{x+2}+a_{x}\,a_{x+1}\,\left(a_{x-1}\,a_{%
 x+1}-2\,a^{2}_{x}\right) $$

to get
$$\twoline{ a_{x-1}\,a^{2}_{2\,x-1}-a_{x}\,\left(2\,a^{2}_{x-1}\,a^{2}_{x+1}-a%
 _{x-1}\,a^{2}_{x}\,a_{x+1}+a^{4}_{x}\right)\,a_{2\,x-1}}{5pt}{+a_{x+1}\,%
 \left(a^{4}_{x-1}\,a^{4}_{x+1}-4\,a^{3}_{x-1}\,a^{2}_{x}\,a^{3}_{x+1%
 }+8\,a^{2}_{x-1}\,a^{4}_{x}\,a^{2}_{x+1}-6\,a_{x-1}\,a^{6}_{x}\,a_{x%
 +1}+2\,a^{8}_{x}\right)=0\rm.} $$

Similarly,
$$ \twoline{a^{3}_{x-1}\,a^{2}_{2\,x}+a_{x}\,\left(2\,a_{x-1}\,a_{x+1}-a^{2}_{%
 x}\right)\,\left(a^{2}_{x-1}\,a^{2}_{x+1}-5\,a_{x-1}\,a^{2}_{x}\,a_{%
 x+1}+a^{4}_{x}\right)\,a_{2\,x}}{5pt}{+a^{3}_{x+1}\,\left(a^{4}_{x-1}\,a^{4%
 }_{x+1}-8\,a^{3}_{x-1}\,a^{2}_{x}\,a^{3}_{x+1}+20\,a^{2}_{x-1}\,a^{4%
 }_{x}\,a^{2}_{x+1}-14\,a_{x-1}\,a^{6}_{x}\,a_{x+1}+3\,a^{8}_{x}%
 \right)=0} $$
and
$$\twoline{ a^{5}_{x-1}\,a^{2}_{2\,x+1}-a_{x}\,\left(8\,a^{4}_{x-1}\,a^{4}_{x+%
 1}-18\,a^{3}_{x-1}\,a^{2}_{x}\,a^{3}_{x+1}+22\,a^{2}_{x-1}\,a^{4}_{x%
 }\,a^{2}_{x+1}-9\,a_{x-1}\,a^{6}_{x}\,a_{x+1}+a^{8}_{x}\right)\,a_{2%
 \,x+1}}{5pt}{+a^{5}_{x+1}\,\left(4\,a^{4}_{x-1}\,a^{4}_{x+1}-12\,a^{3}_{x-1%
 }\,a^{2}_{x}\,a^{3}_{x+1}+20\,a^{2}_{x-1}\,a^{4}_{x}\,a^{2}_{x+1}-16%
 \,a_{x-1}\,a^{6}_{x}\,a_{x+1}+7\,a^{8}_{x}\right)=0\rm.} $$

The square roots plus the size of these expressions probably render them ``too
cumbrous to be of any importance,'' but the even coefficients may pay off
in some finite fields.

``Derivation'' of Conjecture 4a:  By Conjecture 4.5,

\begin{eqnarray*}D_a\pmatrix{k/2-3/2,&k/2-1/2,&k/2+1/2\cr k/2-1/2,&k/2+1/2,&k/2+3/2}&=&
 \left|\matrix{a_{-1}\,a_{k-2}&a_{-2}\,a_{k-1}&a_{-3}\,a_{k}\cr a_{0}\,a%
 _{k-1}&a_{-1}\,a_{k}&a_{-2}\,a_{k+1}\cr a_{1}\,a_{k}&a_{0}\,a_{k+1}&%
 a_{-1}\,a_{k+2}\cr }\right|\;=\;0\cr
\end{eqnarray*}
\begin{eqnarray*}
 &=&2\,\bigl((4\,a_{k-2}\,a_{k}-3\,a^{2}_{k-1})\,a_{k+2}-3\,a_{k-2}\,a^{2}_{k+1}
 +8\,a_{k-1}\,a_{k}\,a_{k+1}-%
 7\,a^{3}_{k}\bigr).\cr
\end{eqnarray*}

\begin{eqnarray*}D_a\pmatrix{ k/2-1/2,&k/2+1/2,&k/2+3/2\cr k/2-3/2,&k/2-1/2,&k/2+1/2\cr}&=&
 \left|\matrix{a_{1}\,a_{k-2}&a_{0}\,a_{k-1}&a_{-1}\,a_{k}\cr a_{2}\,a_{%
 k-1}&a_{1}\,a_{k}&a_{0}\,a_{k+1}\cr a_{3}\,a_{k}&a_{2}\,a_{k+1}&a_{1%
 }\,a_{k+2}\cr }\right|\;=\;0\cr
\noalign{\vskip5pt}
&=&\left(a_{k-2}\,a_{k}-a^{2}_{k-1}\right)\,a_{k+2}-a_{%
 k-2}\,a^{2}_{k+1}
 +3\,a_{k-1}\,a_{k}\,a_{k+1}-2\,a^{3}_{k}.\cr
\end{eqnarray*}

Eliminating $a_{k-2}$,
$$ 2\,a_{k}\,\left(a^{2}_{k-1}\,a^{2}_{k+2}-4\,a_{k-1}\,a_{k}\,a_{k+1%
 }\,a_{k+2}+a^{3}_{k}\,a_{k+2}+a_{k-1}\,a^{3}_{k+1}+a^{2}_{k}\,a^{2}%
 _{k+1}\right)=0,$$
as desired.

{\bf Nonstandard initialization:}  You might wonder how this third order recurrence
can compute a fourth order recurrence with four initial conditions ($a_0,\ldots a_3=1$).
First of all, given the palindrome condition and scaling, there is only one degree
of freedom.  {\it I.e.}, in general, we have

$$a_{-1},a_0,\ldots= p^{2}+p,1,p,p,1,p^{2}+p,p^{2}+p+{{1}\over{p}},
 p^{3}+2\,p^{2}+2\,p+{{1}\over{p^{2}}}+1,\ldots\ .\eqno{\rm(1pp1)}$$

When $p$ is a root of unity, the denominators remain bounded and can be scaled out,
{\it e.g.,}
$$ \ldots,i-1,1,i,i,1,i-1,-1,i-2,2-3\,i,-1,13\,i+3,-16\,i-15,-19\,i-44,\ldots$$
$$ \displaylines{\ldots,i\left(\sqrt{2}+1\right)+1,\sqrt{2},i+1,i+1,\sqrt{2},i\left(\sqrt{2}+1\right)+1,i\,%
 \sqrt{2}+2,i\left(\sqrt{2}+3\right)+\sqrt{2}+1,\cr
 i\left(3\sqrt{2%
 }+7\right)-\sqrt{2}-3,i\left(5\sqrt{2}+12\right)+4\sqrt{2}+2,i%
 \left(15\sqrt{2}+13\right)-22\sqrt{2}-31,i\left(8\sqrt{2}+%
 8\right)-43\sqrt{2}-76,\cr
 i\left(34\sqrt{2}+57\right)-190\sqrt{%
 2}-287,\ldots} $$
$$ \displaylines{\ldots,2\,i\,\sqrt{3},2,i\,\sqrt{3}+1,i\,\sqrt{3}+1,2,2\,i\,%
 \sqrt{3},i\,\sqrt{3}+1,3\,i\,\sqrt{3}-1,\cr
-10,6\,i\,\sqrt{3}-8,-21\,i%
 \,\sqrt{3}-9,35-9\,i\,\sqrt{3},136-66\,i\,\sqrt{3},\ldots}  $$
$$ \displaylines{\ldots, i\left(\sqrt{6}-\sqrt{2}+2\right)+\left(\sqrt{2}+1%
 \right)\sqrt{6}+\sqrt{2},4,i\left(\sqrt{6}-\sqrt{2}\right)+%
 \sqrt{6}+\sqrt{2},i\left(\sqrt{6}-\sqrt{2}\right)+\sqrt{6}+\sqrt{2%
 },4,\cr
 i\left(\sqrt{6}-\sqrt{2}+2\right)+\left(\sqrt{2}+1\right)\,%
 \sqrt{6}+\sqrt{2},2\,i+\left(\sqrt{2}+2\right)\sqrt{6}+2\sqrt{2},\cr
 i\left(2\sqrt{6}+2\right)+\left(3\sqrt{2}+2\right)\sqrt{6}+%
 4\sqrt{2}+4,i\left(\left(4\sqrt{2}+3\right)\sqrt{6}+11\,%
 \sqrt{2}+6\right)+\left(5\sqrt{2}+5\right)\sqrt{6}+15\sqrt{2}+20,\cr
 i\left(\left(5\sqrt{2}+6\right)\sqrt{6}+6\sqrt{2}+10%
 \right)+\left(15\sqrt{2}+18\right)\sqrt{6}+38\sqrt{2}+46,\cr
 i\left(\left(46\sqrt{2}+67\right)\sqrt{6}+107\sqrt{2}+156%
 \right)+\left(54\sqrt{2}+89\right)\sqrt{6}+131\sqrt{2}+212,\cr
 i\left(\left(210\sqrt{2}+311\right)\sqrt{6}+523\sqrt{2}+772%
 \right)+\left(250\sqrt{2}+341\right)\sqrt{6}+635\sqrt{2}+860,\cr
 i\left(\left(963\sqrt{2}+1410\right)\sqrt{6}+2346\sqrt{2}+%
 3434\right)+\left(1383\sqrt{2}+1934\right)\sqrt{6}+3394\sqrt{2%
 }+4754,\ldots\cr} $$
$$ \displaylines{\ldots, i\left(\sqrt{2}\sqrt{5}-\sqrt{2}\right)+\left(i\,%
 \left(\sqrt{5}-1\right)+2\right)\sqrt{\sqrt{5}+5}+\sqrt{2}\sqrt{%
 5}+\sqrt{2},4\sqrt{2},\cr
 2\sqrt{\sqrt{5}+5}+i\sqrt{2}\left(%
 \sqrt{5}-1\right),2\sqrt{\sqrt{5}+5}+i\sqrt{2}\left(\sqrt{5}-1%
 \right),4\sqrt{2},\cr
 i\left(\sqrt{2}\sqrt{5}-\sqrt{2}\right)+%
 \left(i\left(\sqrt{5}-1\right)+2\right)\sqrt{\sqrt{5}+5}+\sqrt{2%
 }\sqrt{5}+\sqrt{2},\cr
 \left(i\left(\sqrt{5}-1\right)+4\right)\,%
 \sqrt{\sqrt{5}+5}+\sqrt{2}\sqrt{5}+\sqrt{2},\cr
 i\left(3\sqrt{2}\,%
 \sqrt{5}-\sqrt{2}\right)+\sqrt{\sqrt{5}+5}\left(\sqrt{5}+i\left(%
 \sqrt{5}-1\right)+3\right)+3\sqrt{2}\sqrt{5}+7\sqrt{2},\cr
 \sqrt{\sqrt{5}+5}\left(i\left(3\sqrt{5}+5\right)+6\sqrt{5}+2%
 \right)+i\left(8\sqrt{2}\sqrt{5}+8\sqrt{2}\right)+9\sqrt{2%
 }\sqrt{5}+13\sqrt{2},\cr
 \sqrt{\sqrt{5}+5}\left(i\left(5\sqrt{%
 5}+5\right)+12\sqrt{5}+20\right)+12\,i\sqrt{2}\sqrt{5}+20\,%
 \sqrt{2}\sqrt{5}+46\sqrt{2},\cr
 \sqrt{\sqrt{5}+5}\left(i\left(37%
 \sqrt{5}+95\right)+47\sqrt{5}+65\right)+i\left(67\sqrt{2}\,%
 \sqrt{5}+187\sqrt{2}\right)+81\sqrt{2}\sqrt{5}+137\sqrt{2},\cr
 \sqrt{\sqrt{5}+5}\left(i\left(168\sqrt{5}+434\right)+162\,%
 \sqrt{5}+322\right)+i\left(332\sqrt{2}\sqrt{5}+796\sqrt{2}%
 \right)+285\sqrt{2}\sqrt{5}+655\sqrt{2},\cr
 \sqrt{\sqrt{5}+5}\,%
 \left(i\left(753\sqrt{5}+1957\right)+883\sqrt{5}+1865\right)+i%
 \left(1499\sqrt{2}\sqrt{5}+3579\sqrt{2}\right)+1644\sqrt{2%
 }\sqrt{5}+3634\sqrt{2},\ldots.\cr}  $$

However, for other $p$ on the unit circle, it is impossible to scale out denominators,
even with a geometric progression.  {\it E.g.},

$$ \displaylines{\ldots, {{39\,i}\over{5}}+\mathchoice {{27}\over{5}}{{27}\over{5}}%
 {27/5}{27/5},5,3\,i+4,3\,i+4,5,{{39\,i}\over{5}}+\mathchoice {{27%
 }\over{5}}{{27}\over{5}}{27/5}{27/5},{{24\,i}\over{5}}+\mathchoice %
 {{47}\over{5}}{{47}\over{5}}{47/5}{47/5},\cr
 {{387\,i}\over{25}}+%
 \mathchoice {{386}\over{25}}{{386}\over{25}}{386/25}{386/25},{{36783%
 \,i}\over{625}}+\mathchoice {{3494}\over{625}}{{3494}\over{625}}{%
 3494/625}{3494/625},\cr
 {{2088\,i}\over{25}}+\mathchoice {{48643}\over{%
 625}}{{48643}\over{625}}{48643/625}{48643/625},{{32642451\,i}\over{%
 78125}}-\mathchoice {{12616857}\over{78125}}{{12616857}\over{78125}}%
 {12616857/78125}{12616857/78125},{{2038374144\,i}\over{1953125}}-%
 \mathchoice {{1297298183}\over{1953125}}{{1297298183}\over{1953125}}%
 {1297298183/1953125}{1297298183/1953125},\cr
 {{9109927539\,i}\over{%
 1953125}}-\mathchoice {{22447118648}\over{9765625}}{{22447118648%
 }\over{9765625}}{22447118648/9765625}{22447118648/9765625},\ldots\cr}  $$

$$ \displaylines{\ldots, {{185\,i}\over{13}}+\mathchoice {{275}\over{13}}{{275%
 }\over{13}}{275/13}{275/13},13,5\,i+12,5\,i+12,13,{{185\,i}\over{13%
 }}+\mathchoice {{275}\over{13}}{{275}\over{13}}{275/13}{275/13},{{%
 120\,i}\over{13}}+\mathchoice {{431}\over{13}}{{431}\over{13}}{431/%
 13}{431/13},\cr
 {{5285\,i}\over{169}}+\mathchoice {{11722}\over{169}}{{%
 11722}\over{169}}{11722/169}{11722/169},{{4966665\,i}\over{28561}}+%
 \mathchoice {{4473518}\over{28561}}{{4473518}\over{28561}}{4473518/%
 28561}{4473518/28561},\cr
 {{468840\,i}\over{2197}}+\mathchoice {{%
 14119307}\over{28561}}{{14119307}\over{28561}}{14119307/28561}{%
 14119307/28561},{{138123747005\,i}\over{62748517}}+\mathchoice {{%
 79921772175}\over{62748517}}{{79921772175}\over{62748517}}{%
 79921772175/62748517}{79921772175/62748517},{{97867476702880\,i%
 }\over{10604499373}}+\mathchoice {{46503584483049}\over{10604499373%
 }}{{46503584483049}\over{10604499373}}{46503584483049/10604499373}{%
 46503584483049/10604499373},\cr
 {{438606726214525\,i}\over{10604499373}}%
 +\mathchoice {{3668262036619888}\over{137858491849}}{{%
 3668262036619888}\over{137858491849}}{3668262036619888/137858491849}%
 {3668262036619888/137858491849},\ldots,}  $$
wherein the powers of $1/5$ and $1/13$ grow quadratically.

Interestingly, (1pp1) permits definition of sequences containing 0.  {\it E.g.}, $p:=-1$
gives a period 10 sequence with values in $\{-1,0,1\}$.  Alternatively, if
$p^3+p^2+1=p,$ the sequence is period 5.  The appropriate generalization of Conjecture 4a
is
$$ a_{k-1}\,a_{k}\,a_{k+1}\,a_{k+2}\,\left(p^{4}+2\,p^{3}+1\right)=%
 \left(a^{2}_{k-1}\,a^{2}_{k+2}+a^{3}_{k}\,a_{k+2}+a_{k-1}\,a^{3}_{k+%
 1}+a^{2}_{k}\,a^{2}_{k+1}\right)\,p^{2}\rm. $$

There appears to be a generalization of this relation for initializations in violation
of the palindrome property.

Perhaps the most important of these is Sloane's A051138

$$A_{-1},A_0,\ldots=-1,0,1,1,-1,-5,-4,29,129,-65,\ldots$$
where
$$A_n=-A_{-n}={A_{n-1}A_{n-3}+A_{n-2}^2\over A_{n-4}}=
 {{A_{n-1}\,a_{n+1}-A_{n-2}\,a_{n+2}}\over{a_{n}}}={{A_{n-1}\,a_{k-%
  1}\,a_{k-n}-A_{n-2}\,a_{k}\,a_{k-n-1}}\over{a_{k-2}\,a_{k-n+1}}}.$$

We can think of $A_n$ as $\sinh$ and $a_n$ as $\cosh$, but actually they're both
theta functions.  Also, $a_n$ is centered at $n=3/2$ instead of $0$.

Solving this last equation for $a_k$ generalizes the Somos4 defining recurrence:

$$a_{k}={{A_{n+1}\,a_{k-1}\,a_{k-n-2}-A_{n+2}\,a_{k-2}\,a_{k-n-1}}\over{A_{n}\,a_{k-n-3}}}\rm.$$

Another such generalization is the ``$k$-tuple speedup'':

$$A^{2}_{k}\,a_{n}\,a_{n+4\,k}= A^{2}_{2\,k}\,a_{n+k}\,a_{n+3\,k}-A_{%
 k}\,A_{3\,k}\,a^{2}_{n+2\,k}\rm. $$

Generalizing both of these is the three-variable relation
$$ a_{n}={{A_{2\,j}\,A_{k+j}\,a_{n-j}\,a_{n-k-2\,j}-A_{j}\,A_{k+2\,j}%
 \,a_{n-2\,j}\,a_{n-k-j}}\over{A_{j}\,A_{k}\,a_{n-k-3\,j}}}. $$

These expressions mixing $A$ and $a$ are somewhat striking because up until now,
all the monomials in a given relation have had the same subscript sum, modulo
$A_{n}=-A_{-n}$ and $a_n=a_{3-n}$.  In particular, these nonconforming identities
can not come directly from $D_s$ type determinant identities, except via the
artifice of multiplying the deficient monomials by $A_k$ and the overweight
monomials by $-A_{-k}$.

$A$ can be eliminated from the speedup identity via the relations
\begin{eqnarray*} {{A_{3\,k}}\over{A_{k}}}&=&{{a_{p}\,a_{p+k+1}\,a_{p+3\,k+1}\,a_{p+4%
 \,k}-a_{p+1}\,a_{p+k}\,a_{p+3\,k}\,a_{p+4\,k+1}}\over{a_{p+k}\,a^{2}%
 _{p+2\,k+1}\,a_{p+3\,k}-a_{p+k+1}\,a^{2}_{p+2\,k}\,a_{p+3\,k+1}}}\rm,\\
 {{A^{2}_{2\,k}}\over{A^{2}_{k}}}&=&{{a_{m}\,a^{2}_{m+2\,k+1}\,a_{m+4%
 \,k}-a_{m+1}\,a^{2}_{m+2\,k}\,a_{m+4\,k+1}}\over{a_{m+k}\,a^{2}_{m+2%
 \,k+1}\,a_{m+3\,k}-a_{m+k+1}\,a^{2}_{m+2\,k}\,a_{m+3\,k+1}}}\rm,\cr
\end{eqnarray*}

for arbitrary $m$ and $p$.  Also,
$$ A^{2}_{k}=a_{k}\,a_{k+3}-a_{k+1}\,a_{k+2} \rm.$$

If we eliminate $A$ between this and the $k$-tuple speedup identity, we
get a polynomial in $a$ with subscript sums which can be brought into
agreement via selective application of $a_n=a_{3-n}$.

Also, $A_n=s_{2n},$ where $s_n=$ Sloane's A006769:
$$s_{-1},s_0,\ldots=-1,0,1,1,-1,1,2,-1,-3,-5,7,-4,-23,29,\ldots,$$

and $$s_n = -s_{-n} = {s_{n-1}s_{n-3}+s_{n-2}^2\over s_{n-4}}\rm,$$
the same recurrence as $A_n$.  Perhaps surprisingly,
$$ s_{2\,n+1}=(-1)^na_{n+2}\rm. $$

That $A_{1/2},A_{3/2},\ldots$ can be integers suggests that $a_{1/2},
a_{3/2},\ldots$ could be, too.  Substituting half-integers into the
$\vartheta$ expression below yields nonintegers, but it is likely that
there are alternative analytic expressions for $a_n$ which disagree for
nonintegers.

Curiously, $A_{2n}$ does not obey the Somos4 recurrence.

Note that $a$ is even easier than $A$ to eliminate from the mixed recurrences,
since they hold for $a=A$!  {\it I.e.,}

$$ A_{j}\,A_{k}\,A_{n}\,A_{n-k-3\,j}=A_{2\,j}\,A_{k+j}\,A_{n-j}\,A_{n%
 -k-2\,j}-A_{j}\,A_{k+2\,j}\,A_{n-2\,j}\,A_{n-k-j}. $$

With the relation $A_{-n}=-A_n$ along with linear changes of variable, this
can be rewritten
$$ A_{-j}\,A_{j-k}\,A_{j-n}\,A_{n+k+j}+A_{2\,j}\,A_{-k}\,A_{-n}\,A_{n%
 +k}=A_{-j}\,A_{k+j}\,A_{-n-k+j}\,A_{n+j}, $$

so that each term's subscript sum is $2j$.  We might thus expect an equivalent
3 by 3 determinant {\it a la} Conjecture 4.  The most general case gives a
six variable relation with 24 terms of degree 6.  The only apparent way to
reduce to degree four is to specialize two of the variables to create terms
of absolute value 1, {\it i.e.} $A_{\pm1},A_{\pm2},$ or $A_{\pm3}$.  But this
will introduce small integer offsets among the remaining subscripts, a feature
notably absent from our trivariate relation.  So for $A_n$, at least, determinants
may not tell the whole story.

Likewise for $\vartheta_1$:  the trivariate relation empirically holds if
we replace $A_n$ by $\vartheta_1(n,q)$, for arbitrary complex $n$ and any
fixed q within the unit circle.  That it fails for the other $\vartheta$s
suggests the existence of a four or more variable generalization.  Indeed,
by analogy with (4vars),
$$ A_{k-i}\,a_{k+i}\,A_{j-n}\,a_{n+j}=A_{j-i}\,a_{j+i}\,A_{k-n}\,a_{n%
 +k}+A_{k-j}\,a_{k+j}\,A_{i-n}\,a_{n+i}, $$
also holding with $A$ in place of $a$.  So maybe (4vars) type determinants {\it do}
tell the whole story.

If $s_n=-s_{-n},\ s_1=1,$ (as with $s_n:=A_n$,) then

\begin{eqnarray*}D_s\pmatrix{y,&0,&1\cr n,&0,&1}&=&
\left|\matrix{-s_{n-y}\,s_{y+n}&s^{2}_{y}&s_{y-1}\,s_{y+1}\cr -s^{2}_{n%
 }&0&-1\cr -s_{n-1}\,s_{n+1}&1&0\cr }\right|\cr
\noalign{\vskip5pt}
&=& -s_{n-y}\,s_{n+y}-s^{2}_{n}\,s _{y-1}\,s_{y+1}+s_{n-1}\,s_{n+1}\,s^{2}_y\;\;=\;\;0
   \;\;\;\;\;\;\;\hspace{1in} \rm(EDS)\cr
\end{eqnarray*}

is equivalent to $s$ being an elliptic divisibility sequence.  Integer divisibility
sequences merely require $d|n\Rightarrow s_d|s_n$, but the divisibility sequences
discussed in this report appear to satisfy the stronger relation
$(s_x,s_y)=|s_{(x,y)}|$, even when they disobey the addition formula.
This may be what is meant by ``strong divisibility sequence''.

The EDS upside is this nice addition formula.

Fomin and Zelevinsky have shown that Somos4, $\ldots$, Somos7 are Laurent
polynomials (rational functions with monomial denominators) in their initial values.

{\bf Somos4oid polynomials:}
We can get true polynomials from the ``odd'' ($s_{-n}=-s_n$) sequences with the
initialization $-1, 0, 1, 1, -1, x$, where $x$ is unconstrained by the Somos4
recurrence, which gives $0/0$.  At greater length,

$$\twoline{ s_n= -1,0,1,1,-1,x,x+1,x^{2}-x-1,-x^{3}-x-1,-3\,x^{2}-2\,x,x^{5%
 }-x^{4}+3\,x^{2}+3\,x+1,}{5pt}{-x^{6}-2\,x^{4}-5\,x^{3}+3\,x+1,-x^{7}+2\,x%
 ^{6}-3\,x^{5}-9\,x^{4}-5\,x^{3}-3\,x^{2}-3\,x-1,\ldots.} $$

$A_n$ is the case $x=-5$ and $A_{n/2}$ is the case $x=1$.  {\it I. e.},
$s_n(-5)=s_{2n}(1)$.  Empirically, this is a
strong (redundant?) elliptic polynomial(!) division sequence for all
$x$.  If indeed the divisibility property holds for both integers and polynomials,
then the values assumed by the polynomials $s_k(x)/s_{(k,n)}(x)$
and $s_n(x)/s_{(k,n)}(x)$ are relatively prime for every integer x.

It seems that the Chebychev polynomials $U_{n-1}(y):=\sin(n\arccos y)
/\sin(\arccos y)$ behave similarly.  {\it E.g.}, for integers $k$ and $n$,
$\sin(kn\arccos y)/\sin(n\arccos y)$ is a polynomial in $y$, but of degree
only $(k-1)n$.

Here are the polynomial factorizations of $s_n(x)$ through $n=18$.

$$ \matrix{n&s_n\cr\cr
-1&-1\cr 0&0\cr 1&1\cr 2&1\cr 3&-1\cr 4&x\cr 5&x+1\cr 6&x%
 ^{2}-x-1\cr 7&-\left(x^{3}+x+1\right)\cr 8&-x\,\left(3\,x+2\right)%
 \cr 9&x^{5}-x^{4}+3\,x^{2}+3\,x+1\cr 10&-\left(x+1\right)\,\left(x^{%
 5}-x^{4}+3\,x^{3}+2\,x^{2}-2\,x-1\right)\cr 11&-\left(x^{7}-2\,x^{6}%
 +3\,x^{5}+9\,x^{4}+5\,x^{3}+3\,x^{2}+3\,x+1\right)\cr 12&-x\,\left(x%
 ^{2}-x-1\right)\,\left(x^{6}+2\,x^{4}+5\,x^{3}+9\,x^{2}+9\,x+3%
 \right)\cr 13&x^{10}+x^{9}+13\,x^{7}+12\,x^{6}-6\,x^{5}+16\,x^{3}+15%
 \,x^{2}+6\,x+1\cr 14&-\left(x^{3}+x+1\right)\,\left(x^{9}-3\,x^{8}+x%
 ^{7}+6\,x^{6}-13\,x^{5}-30\,x^{4}-15\,x^{3}+4\,x^{2}+5\,x+1\right)%
 \cr
 15&-\left(x+1\right)\,\left(x^{13}-2\,x^{12}+5\,x^{11}-2\,x^{10}%
 +5\,x^{9}+8\,x^{8}\right.\qquad\qquad\qquad\cr
&\qquad\qquad\qquad\left.-19\,x^{7}+12\,x^{6}+63\,x^{5}+50\,x^{4}+20\,x^{3}%
 +10\,x^{2}+5\,x+1\right)\cr 16&x\,\left(3\,x+2\right)\,\left(2\,x^{%
 12}-3\,x^{11}+6\,x^{10}+14\,x^{9}-2\,x^{8}\right.\qquad\qquad\qquad\cr
 &\qquad\qquad\qquad\left.+3\,x^{7}+23\,x^{6}+18\,x%
 ^{5}-6\,x^{4}-27\,x^{3}-27\,x^{2}-12\,x-2\right)\cr 17&x^{18}-3\,x^{%
 17}+4\,x^{16}+6\,x^{15}-9\,x^{14}+5\,x^{13}+56\,x^{12}+69\,x^{11}+%
 105\,x^{10}+311\,x^{9}\cr
&\qquad+429\,x^{8}+211\,x^{7}-2\,x^{6}+45\,x^{5}+135%
 \,x^{4}+110\,x^{3}+45\,x^{2}+10\,x+1\cr 18&\left(x^{2}-x-1\right)\,%
 \left(x^{5}-x^{4}+3\,x^{2}+3\,x+1\right)\,\left(x^{13}+x^{12}+7\,x^{%
 11}+19\,x^{10}+25\,x^{9}\right.\qquad\qquad\cr
&\qquad\qquad\qquad\qquad\left.+78\,x^{8}+133\,x^{7}+108\,x^{6}+79\,x^{5}+%
 65\,x^{4}+24\,x^{3}-6\,x^{2}-6\,x-1\right)\cr } $$

The degrees of the polynomials, starting with $n=1$, go
$$0,0,0,1,1,2,3,2,5,6,7,9,10,12,14,14,18,20,22,25,27,30,33,34,39,42,45,49,52,56,60,62,68,72,76,81,\ldots,$$

which is eight interlaced quadratic progressions:

$$\deg s_{8q+r}=(4q+r)q+[-2,0,0,0,1,1,2,3]_r,\qquad 0\le r\le7,$$
which can be written
$$\deg s_n = {\sqrt2\over4}(\sin {n\pi\over2})(\sin {n\pi\over4})
 -(\cos {n\pi\over2})({3\over8}+\cos {n\pi\over4})
 +{1\over32}(2n^2-5(-1)^n-5).$$

It appears that $n$ prime $\Rightarrow$ $s_n$ irreducible.  The polynomials appear
to be monic except for $s_{8n}$, whose leading coefficients appear to be $(-)^n3n$.

It appears that all the polynomials $s_n(x)$ have a root close to $x=\omega\approx
-0.669499628215, \; 
\omega^{4}+3\,\omega^{3}-5\,\omega^{2}+21\,\omega+17=0$, with proximity rapidly increasing with $n$.

Besides the EDS condition, we retain the $\vartheta_1$ three-variable identity
$$ s_{2\,j}\,s_{k}\,s_{n}\,s_{n+k}=s_{j}\,s_{k-j}\,s_{n-j}\,s_{n+k+j}%
 +s_{j}\,s_{k+j}\,s_{n+j}\,s_{n+k-j} $$

This can be subscript-balanced as
$$ s_{2\,j}\,s_{k+j}\,s_{-n-k-j}\,s_{n}=s_{-j}\,s_{-k}\,s_{j-n}\,s_{n%
 +k+2\,j}-s_{-j}\,s_{k+2\,j}\,s_{-n-k}\,s_{n+j}, $$
but its asymmetry and failure to subsume the EDS condition suggest that we're
missing a nice, {\it four}-variable relation.  Sure enough, by analogy with (4vars),
$$ s_{k-i}\,s_{k+i}\,s_{j-n}\,s_{n+j}=s_{j-i}\,s_{j+i}\,s_{k-n}\,s_{n%
 +k}+s_{k-j}\,s_{k+j}\,s_{i-n}\,s_{n+i} $$

withstands empirical testing.

This identity specializes to a particularly attractive doubling formula:
\begin{eqnarray*}s_{2\,n-1}&=&s^{3}_{n-1}\,s_{n+1}-s_{n-2}\,s^{3}_{n}\\
 s_{2\,n}&=&\left(s^{2}_{n-1}\,s_{n+2}-s_{n-2}\,s^{2}_{n+1}\right)s_{n}\,.\cr \end{eqnarray*}
Given the four consecutive values $s_{n-2},\ldots,s_{n+1}$, extend them to $s_{n+3}$
stepping the recurrence twice.  Then use the doubling formula to get the four values
$s_{2n-1},\ldots,s_{2n+2}$.  Etc.

We also have
$$ s_{b-a}\,s_{b+a}=s^{2}_{a}\,s_{b-1}\,s_{b+1}-s_{a-1}\,s_{a+1}\,s^{2}_{b}. $$

And we have the ntuple speedup relation
$$ s^{2}_{k}\,s_{n}\,s_{n+4\,k}=s^{2}_{2\,k}\,s_{n+k}\,s_{n+3\,k}-s_{%
 k}\,s_{3\,k}\,s^{2}_{n+2\,k}.$$

This provides an alternative doubling process:  Given four values
$s_k, s_{2k},s_{3k},s_{4k},$ start $n$ at $k$ and generate
$s_{5k},s_{6k},s_{7k},s_{8k}$.  Discard the odd multiples, and we have
doubled $k$ and are free to iterate.

For $x=0$, $s_n$ has period 8:
\begin{eqnarray*}s_{4q+r}(0)&=&-1,0,1,1,-1,0,1,-1,-1,0,1,1,\ldots\\
&=&[0,1,(-1)^q,-1]_r,\qquad\qquad0\le r\le3,\\
&=&bi^{-n^2/4}\vartheta_1(n\pi/4,Q),\qquad n=4q+r,\ b=.2653512762412i+.4652895036579,\\
   && \qquad \qquad \qquad \qquad \qquad \qquad Q=.7359196601139i+.3006597280279.\cr \end{eqnarray*}

Of course, a much simpler expression is
$$ s_{n}\left(0\right)=\sin \left({{\pi\,n}\over{2}}\right)-\sin %
 \left({{\pi\,n}\over{4}}\right)\,\cos \left({{\pi\,n}\over{2}}%
 \right). $$

For $x=-1$ the period is 5:
\begin{eqnarray*}s_{5q+r}(-1)&=&-1,0,1,1,-1,-1,0,1,1,\ldots\\
&=&[0,1,1,-1,-1]_r,\qquad\qquad0\le r\le4,\\
&=&{b\over\root4\of Q}\vartheta_1(2n\pi/5,Q)
 \qquad n=5q+r,\ b=0.6155370356317,\ Q=.4856907848670i.\cr
\end{eqnarray*}

For $x=-2/3$, we get eight interlaced progressions:
\begin{eqnarray*}s_{8q+r}(-2/3)&=&-1,0,1,1,-1,-2/3,1/3,1/3^2,-1/3^3,0,1/3^5,-1/3^6,-1/3^7,
2/3^9,\ldots\\
&=&{[0,3^{3q+1},(-)^q3^{2q+1},-3^{q+1},(-)^{q+1}2,3^{1-q},(-)^q3^{-2q-1},-3^{1-3q}]_r
\over3^{(2q+1)^2}},\qquad\qquad0\le r\le7,\\
 &=&bu^{n^2}\vartheta_1(ny,Q),\qquad n=8q+r,\ b= -2.010659335767i,\ u=0.8509811643954i,\\
 & &\qquad\qquad\qquad\qquad y=\pi/2-.7416161288587i,\ Q = .0026507066057.\cr
\end{eqnarray*}

Note that $y$ is not $\pi/8$ nor even real, so where do the periodic 0s come from?
And $|u|$ is not $3^{-1/16}$, in fact, even its square root is too small.  So
where does the $3^{-n^2/16}$ ``growth'' rate come from?  The answer, as usual,
is clear after Jacobi's imaginary transformation:


$$ s_{n}(-2/3)={{2\,e^{{{9\,i\,\pi\,n^{2}}/{16}}}\,\vartheta_{1}\left(%
 {{\pi\,n}\over{8}}\right)}\over{\vartheta_{2}\,3^{%
 {{n^{2}}/{16}}}}}, $$
with $q$ satisfying
$$ \prod_{i=1}^{\infty }{\left(1+q^{2\,i}\right)\,\left(1+q^{4\,i-2}%
 \right)}={{1+i}\over{\root 4 \of{3}}}, $$
{\it e.g.},
$$ q\approx .5913080374704560258502159338438+%
 .4423170132359810537349781037012\,i .$$

We thus answer both questions, and reduce four mysterious parameters to one.
Or rather, two, since
$$ q\approx.7241830710727415040344246937315\,i+%
 .5068861260317593704061905537186 $$
also satisfies the infinite product constraint, but produces the mysterious
sequence
$$ \ldots,-\sqrt{3\,i},0, \sqrt{3\,i},1,-{{\sqrt{3\,i}}\over{3}},-\mathchoice {{2%
 }\over{3}}{{2}\over{3}}{2/3}{2/3},{{\sqrt{3\,i}}\over{9}},%
 \mathchoice {{1}\over{9}}{{1}\over{9}}{1/9}{1/9},-{{\sqrt{3\,i}%
 }\over{27}},0,{{\sqrt{3\,i}}\over{243}},-\mathchoice {{1}\over{729}}%
 {{1}\over{729}}{1/729}{1/729},-{{\sqrt{3\,i}}\over{6561}},\ldots\ !  $$

For $x=$ the golden ratio, we get six interlaced progressions:
\begin{eqnarray*}s_{6q+r}(\phi)&=&-1,0,1,1,-1,\phi,\phi^2,0,-\phi^4,-\phi^5,\phi^6,-\phi^8,
-\phi^{10},0,\ldots\\
&=&(-)^q\phi^{6{q+1\choose2}}[0,\phi^{-2q},\phi^{-q},-1,\phi^{q+1},\phi^{2q+2}]_r,
\qquad\qquad0\le r\le5.\cr
\end{eqnarray*}

Likewise for the conjugate:
\begin{eqnarray*}s_{6q+r}(-1/\phi)&=&-1,0,1,1,-1,-1/\phi,\phi^{-2},0,-\phi^{-4},
\phi^{-5},\phi^{-6},-\phi^{-8},-\phi^{-10},0,\ldots\\
&= & {[0,\phi^{2q},(-\phi)^q,-1,(-\phi)^{-q-1},\phi^{-2q-2}]_r\over(-)^q\phi^{6{q+1\choose2}}},
\qquad\qquad0\le r\le5,\\
 &=&bu^{n^2}\vartheta_1(ny,Q),\qquad n=6q+r,\ b= -1.403592671340i,\ u=0.8476983265649i,\\
 & &\qquad\qquad\qquad\qquad y=\pi/2-.7507768082213i,\ Q = .0110573396552.\cr
\end{eqnarray*}

This $\vartheta$ expression is close to $s_n(-2/3)$ because $-1/\phi=-.618$ is close to
$-2/3$.  

It is probable that $s_n(\alpha)$ comes out in $k$ such interlaced progressions
when $s_k(\alpha)=0$.  {\it E.g.}, when $\alpha^3+\alpha+1=0$, we appear to get seven
interlaced progressions scaled by $\alpha^{(k+1)(k+1/7)}$.

An algebraic $x$ for which $s_n(x)$ has an elementary closed form is $x=-w^3-2w^2$, where $w^4-w=1$,
{\it i.e.} 
$ x^{4}+3\,x^{3}-5\,x^{2}+21\,x+17=0,\; w=-(3\,x^{3}-26\,x^{2}-89\,x+158)/283 $:
\begin{eqnarray*}s_{n}&=& -1,0,1,1,-1,x,x+1,x^{2}-x-1,-x^{3}-x-1,-x\,\left(3\,%
 x+2\right),17\,x^{3}-38\,x^{2}+70\,x+69,\\
  & &\qquad\qquad79\,x^{3}-126\,x^{2}+288\,x+ 273,
  -526\,x^{3}+925\,x^{2}-1838\,x-1803,\ldots\\
  &=& -1,0,1,1,-1,-w^{3}-2\,w^{2},-w^{3}-2\,w^{2}+1,2\,w^{3}+7%
 \,w^{2}+8\,w+3,22\,w^{3}+18\,w^{2}+25\,w+17,\\
  & &\qquad-w^{3}-11\,w^{2}-24\,w-%
 12,-465\,w^{3}-602\,w^{2}-729\,w-389,\\
  & &\qquad-2073\,w^{3}-2470\,w^{2}-2983\,%
 w-1653,13809\,w^{3}+16717\,w^{2}+20550\,w+11365,\ldots\\
 &=& {{\displaystyle\sin \left(n\,\arccos {w^{-3/2}\over2}\right)\,w^{{{\left(n-1%
 \right)\,\left(n+1\right)}\over{2}}}}\over{\displaystyle\sqrt{1-{w^{-3%
 }\over{4}}}}}\;\;=\;\;U_{n-1}\left({w^{-3/2}\over{2}}\right)\,w^{{{\left(n-1\right)\,\left(n+1%
 \right)}\over{2}}},\\ \end{eqnarray*}
where $U_n$ is the Chebychev polynomial, second kind.  The degeneration of the
$\vartheta$ corresponds to the vanishing of $q$.  Note that one of the roots
$x\approx-0.669499628215$, which is numerically close to $-2/3$, which at least explains
the unusually small value of $q$ in the otherwise puzzling $\vartheta$ expressions
for $s_n(-2/3)$ and $s_n(-1/\phi)$.

This $\vartheta$-free expression for $s_n(-w^3-2w^2)$ affords elementary expansions of $b,u,y,$ and $q$ about $\epsilon=0$ in

$$s_n(\epsilon-w^3-2w^2)
={b\over2\root8\of q}\,u^{n^{2}}\,\vartheta_{1}\left(n\,y,\sqrt q\right)+O(\epsilon^3), $$
namely

\begin{eqnarray*} b&=&{\displaystyle1+{{37\,\epsilon}\over{757\,w^{3}+450\,w^{2}+87\,w-661}}
+{{858967 \,\epsilon^{2}}\over{3241613\,w^{3}+34180558\,w^{2}+150075980\,w+%
 30712993}}+\cdots\over{
    \displaystyle\sqrt{w-{{1}\over{4\,w^{2}}}}}}\\
   &=&{\displaystyle1-{{37\,\epsilon}\over{31\,x^{3}+86\,x^{2}+192\,x+875}}-{{858967\,%
 \epsilon^{2}}\over{807948\,x^{3}-16234660\,x^{2}-29959955\,x+48768711}}+\cdots\over{
    \displaystyle\sqrt{w-{{1}\over{4\,w^{2}}}}}},\\
 u&=&\sqrt{w}\,\left(1-{{\epsilon}\over{16\,w^{3}+32\,w^{2}+2}}-{{2556371\,\epsilon^{2}%
 }\over{1728448\,w^{3}+518695360\,w^{2}+706843904\,w+458996344}}+\cdots\right)\\
  &=&\sqrt{w}\,\left(1+{{\epsilon}\over{16\,x-2}}-{{2556371\,\epsilon^{2}}\over{7072000%
 \,x^{3}+110455488\,x^{2}-39784960\,x+144470392}}+\cdots\right),\\
 y&=&\arccos \left({{\displaystyle1-{{13\,\epsilon}\over{308\,w^{3}+200\,w^{2}+400\,w+258}}-yy+\cdots%
 }\over{2\,w^{\mathchoice {{3}\over{2}}{{3}\over{2}}{3/2}{3/2}}}}%
 \right)\\
  &=&\arccos \left({{\displaystyle1+{{13\,\epsilon}\over{16\,x^{3}-28\,x+30}}+{{290035289\,
\epsilon^{2}}\over{1143740480\,x^{3}-236851392\,x^{2}+7982806560\,x+9574893592}}+\cdots%
 }\over{2\,w^{\mathchoice {{3}\over{2}}{{3}\over{2}}{3/2}{3/2}}}}%
 \right)\!\!,\\
 yy & = & {{
 290035289\,\epsilon^{2}}\over{32238208032\,w^{3}+35449717760\,w^{2}+%
 30488323136\,w+11959840616}},\\
 q&=&{{\epsilon\,\left(\displaystyle1-{{4339\,\epsilon}\over{1733\,w^{3}-55408\,w^{2}+2992\,w
+36867}}+qq+\cdots\right)}\over{234-77\,w+91\,w^{2}-170\,w^{3}}}\\
  &=&{{\epsilon\,\left(\displaystyle1+{{4339\,\epsilon}\over{1696\,x^{3}+4926\,x^{2}-28957\,x-26043}}-%
 qqq +\cdots\right)}\over13\,x^{3}+31\,x^{2}-72\,x+344},\\
 qq & = & {{477927637\,\epsilon^{2}}\over{15857700151\,w^{3}-1431606991\,w^{2}%
 -2016123508\,w+7704208617}},\\
 qqq &=&  {{477927637\,\epsilon^{2}}\over{915645540\,x^{3}+3113407751\,x^{2}-%
 257448438\,x-3676219901}}.\\
\end{eqnarray*}
What is it with $283=566/2=\sqrt{80089}=\root3\of{22665187}$?  Answer:  $-283={\rm
discriminant}(w^4-w-1)=\root3\of{{\rm discriminant}(x^{4}+3\,x^{3}-5\,x^{2}+21\,x+17)}$.

Using the negative root $w\approx-0.72449195900052$, these expansions through $\epsilon^4$ go
\begin{eqnarray*} b&\approx&-0.9125730603509\,i\,(1-0.04769878803144\,\epsilon-0.01401078082448\,\epsilon^{2}\\
   & &\qquad\qquad-0.0050002222638\,\epsilon^{3}-7.44689129426338\cdot 10^{-4}\,%
 \epsilon^{4}+\cdots),\\
 u&\approx&0.85117093406702\,i\,(1-0.07866586437604\,\epsilon-0.01169998061242\,%
 \epsilon^{2}\\
   &  &\qquad\qquad-0.00163136901344\,\epsilon^{3}+%
             1.1997262491263\cdot 10^{-4}\,\epsilon^{4}+\cdots),\\
 y&\approx&\arccos (0.81081103497608\,i\,(1+%
 0.29582732047992\,\epsilon+0.07670823278512\,\epsilon^{2}\\
  &  &\qquad\qquad+0.02128110984092\,\epsilon^{3}+0.00440366306328\,\epsilon^{4}+\cdots)),\\
 q&\approx&0.00248633800734\epsilon(1-0.87527986918762\,\epsilon
+0.20054769983458\,\epsilon^{2}+0.0025905870597\,\epsilon^{3}+\cdots),\\ \end{eqnarray*}

suggesting a fairly commodious radius of convergence.

The positive root $ w=1.22074408460576,
 s_n(\epsilon-w^3-2w^2)=-1,0,1,1,-1,\epsilon-4.7996,\ldots,$ has expansions
\begin{eqnarray*} b&=&.97451669219348\,(1+.02478371881502\,\epsilon+.00317253690858\,\epsilon^2\\
   & &\qquad\qquad+4.97976424690216\cdot 10^{-4}\,\epsilon^{3}+%
   8.74003233413136\cdot 10^{-5}\,\epsilon^{4}+\cdots ),\\
u&=&1.10487288165008\,(1-.01269137385174\,\epsilon-.00121848904668\,\epsilon^{2}\\
   & &\qquad\qquad- 1.63711992744874\cdot 10^{-4}\,\epsilon^{3}-%
 2.5512354068948\cdot 10^{-5}\,\epsilon^{4}+\cdots ),\\
y&=&\arccos \left(.37070894172584\,(1-%
 .00810147538142\,\epsilon-.0018053537704\,\epsilon^{2}\right.\\
  & &\qquad\qquad\left.-3.84723167907628\cdot 10^{-4}\,\epsilon^{3}-%
 8.0356533846957\cdot 10^{-5}\,\epsilon^{4}+\cdots)\right),\\
q&=&-.02972037154846\,\epsilon\,(1+.1115487179065%
 \,\epsilon+.01495508864096\,\epsilon^{2}\\
 & &\qquad\qquad+.00221611060456\,\epsilon ^{3}+\cdots).\\ \end{eqnarray*}


Another interesting Somos4 (apparently (strong) polynomial (non-E)DS) is
\begin{eqnarray*} s_{0},s_{1},\ldots & = & 0,1,i,1,x,i\,\left(x-i\right)%
 ,-i\,\left(x^{2}+i\,x+1\right),-i\,\left(x^{3}-x+i\right),-i\,x\,%
 \left(3\,x-2\,i\right),i\,\\
 & & \left(x^{5}+i\,x^{4}+3\,i\,x^{2}+3\,x-i%
 \right),
  i\,\left(x-i\right)\,\left(x^{5}+i\,x^{4}-3\,x^{3}+2\,i\,x^{%
 2}-2\,x+i\right),\\
 & & -i\,\left(x^{7}+2\,i\,x^{6}-3\,x^{5}+9\,i\,x^{4}+5%
 \,x^{3}-3\,i\,x^{2}-3\,x+i\right),\\
 & & -x\,\left(x^{2}+i\,x+1\right)\,%
 \left(x^{6}-2\,x^{4}+5\,i\,x^{3}+9\,x^{2}-9\,i\,x-3\right),
\end{eqnarray*}
which gives us Gaussian integers, among other things.  As with the previous Somos4
polynomial sequence, there is likely a value of $x$ for which the $\vartheta_1$
degenerates to a Chebychev, and consequently another set of elementary expansions of the
$\vartheta$ parameters about this $x$.  But foo, these polynomials are essentially
identical to those generated by the $-1,0,1,1,-1,x,x+1,\ldots$ sequence.

{\bf Corollary 4:}  the determinant
$$
 D_a\pmatrix{s,&t,&u,&v\cr w,&x,&y,&z\cr}=0,$$
where 
$s,t,u,v,w,x,y,$ and $z$ are arbitrary integers.  Proof:  Dodgson's rule, provided the
central 2 by 2 doesn't vanish.

{\bf Expression as $\vartheta$:}  Email from Noam Elkies to sci.math suggests the
relation
$$ a_{n}=b\,u^{\left(n-\mathchoice {{3}\over{2}}{{3}\over{2}}{3/2}{3/%
 2}\right)^{2}}\,\sum_{k=-\infty }^{\infty }{q^{k^{2}}\,z^{k\,\left(n%
 -\mathchoice {{3}\over{2}}{{3}\over{2}}{3/2}{3/2}\right)}}
   =b\,u^{\left(n-\mathchoice {{3}\over{2}}{{3}\over{2}}{3/2}{3/2}%
 \right)^{2}}\,\vartheta_{3}\left(i\,\left(n-\mathchoice {{3}\over{%
 2}}{{3}\over{2}}{3/2}{3/2}\right)\,{{\log z}\over{2}},q\right) \rm.$$
Using $n\in\{2,3,4,5\}$ to numerically approximate $b,u,q,$ and $z$,
$$ \left\{ b=1.01943271913292,u=0.63853138366726,z=0.05462469648874,q%
 =0.02157360406362\right\} \rm. $$
These constants do not appear to be in Plouffe's collection.
Plugging in $-3,-2,\ldots,10,11$ for $n$ gives
$$ \twoline{\{ 6.99999999999998,3.0,2.0,1.0,1.0,1.0,1.0,2.0,3.0,%
 6.99999999999998,23.0,58.9999999999998,}{4pt}{313.999999999998,%
 1528.99999999998,8208.99999999994\}\rm,}  $$
in good agreement.  Due to the pleasantly small value of $q$, we even get
$1123424582770.98$ for $a_{18}=1123424582771$.  In fact, the only terms affecting
this double precision result were $-9\le n\le-2$, making the series a competitive
alternative numerical method.

With Jacobi's imaginary transformation, we get an even nicer, entirely real
expression:

$$ a_n={{b\,u^{\left(n-\mathchoice {{3}\over{2}}{{3}\over{2}}{3/2}{3/2}%
 \right)^{2}}\,\vartheta_{4}\left(\left(n-\mathchoice {{3}\over{2}}{{%
 3}\over{2}}{3/2}{3/2}\right)\,y,q\right)}} $$
with
$$\twoline{y\approx1.9511889024071,\ q\approx.07632928490026,}{4pt}{
 b={{\vartheta_{4}\left({{3\,y}\over{2}},q\right)^{1/8}}\over{%
 \vartheta_{4}\left({{y}\over{2}},q\right)^{\mathchoice {{9}\over{8}}%
 {{9}\over{8}}{9/8}{9/8}}}}\approx.92252487906093,\ 
 u=\sqrt{{{\vartheta_{4}\left({{y}\over{2}},q\right)}\over{%
 \vartheta_{4}\left({{3\,y}\over{2}},q\right)}}}\approx1.10763024250632.}  $$

$A_n$ not surprisingly comes out as a $\vartheta_{1}$, also with a fairly small $q$:
\begin{eqnarray*} A_{n}&=&-2\,b\,u^{n^{2}}\,\sum_{k=0}^{\infty }{\left(-1\right)^{k}\,q^{%
 \left(2\,k+1\right)^{2}}\,\sin \left(\left(2\,k+1\right)\,n\,y \right)}\cr
 &=&{{-b}}\,u^{n^{2}}\,\vartheta_{1} \left(n\,y,q^{4}\right),\cr
\end{eqnarray*}
where
$$y=1.9511889024071,\ q=0.52562110924304,\ b=.92252487906093,\ u=-1.10763024250632.$$
Note that this $q$ is raised to the fourth power in the $\vartheta_{1}$,
so that, modulo an alternating sign, these parameters are identical to those
in the $\vartheta_4$ formula for $a_n$.
I found these parameters enormously tough to compute (prior to Jacobi-transforming the
$a_n$ expression), which may explain their
 absence from Elkies's email.  Then again I flunked numerical analysis.  Of course, now
that we have it, the Chebychev expansion is also valid at $w=1.22074408460576, x=w^3-2w^2=
 - 4.79960475359606,$ which is close enough to $-5$ to provide an excellent first
approximation.  To a different solution, however!  (Negated $b$ and $u$, $\pi$-complement
 of $y$.)

Testing the non-Chebychev expression:

$0 = 0.0d0, 1 = 1.0d0, 1 = 1.0d0,  - 1 =  - 1.0d0, 
 - 5 =  - 4.99999999999998d0,  - 4 =  - 3.99999999999998d0, 
29 = 28.9999999999998d0, 129 = 128.999999999998d0, 
 - 65 =  - 65.0000000000002d0,  - 3689 =  - 3688.99999999996d0\ .$

Recalling that $A_n=s_n(-5)=s_{2n}(1)$, where $s_n(x)$ is the EDS polynomial
sequence satisfying Somos4, we sought a $\vartheta_1$ expression for $s_n(1)$
to see if $s_{2n}(1)$ gives the same expression as the $\vartheta_1$ for
$A_n$.  In fact, we found (with much difficulty) ten $\vartheta_1$
expressions that agree with $s_n(1)$ for integer $n$.  For noninteger or nonreal
$n$, symmetry suggests that there are as many as sixteen different functions.
(And not one of them coincides, for $n\leftarrow 2n$, with our $A_n$ expression.)
The sixteen seem to divide into two classes of eight.  Within each class,
their values at $n=1/2, 3/2,\ldots$ agree modulo conjugation and multiplication
by some integer power of $i$.

However, there appear to be {\it more} than sixteen $(y,q)$ pairs producing $s_n(1)$.
{\it I.e.}, there are multiple ways to express $s_n(1)$ as $b\,u^{n^2}\,
 \vartheta_1(ny,q)$ that agree even for complex $n$!  In particular,
$$ b=-{{\vartheta_{1}\left(2\,y,q\right)^{3}}\over{\vartheta%
 _{1}\left(y,q\right)^{3}\,\vartheta_{1}\left(3\,y,q\right)}},\quad u=-{{%
 \vartheta_{1}\left(y,q\right)^{2}\,\vartheta_{1}\left(3\,y,q\right)%
 }\over{\vartheta_{1}\left(2\,y,q\right)^{3}}},  $$

and

$$ s_n(1)=-{{\vartheta_{1}\left(2\,y,q\right)^{3}\,\vartheta_{1}%
 \left(n\,y,q\right)}\over{\vartheta_{1}\left(y,q\right)^{3}\,%
 \vartheta_{1}\left(3\,y,q\right)}}\,\left(-{{\vartheta_{1}%
 \left(y,q\right)^{2}\,\vartheta_{1}\left(3\,y,q\right)}\over{%
 \vartheta_{1}\left(2\,y,q\right)^{3}}}\right)^{n^{2}} $$

seems to be exactly the same function of complex $n$ for
$$  y\approx0.49235539271999-0.74875275029651\,i,\quad q\approx%
 0.69018582634555\,i+0.54229640598463$$
as for
$$ y\approx 1.11453161008963\,i+0.62943384983216,\quad q\approx0.63418111840451\,i-%
 0.43035475675355. $$
(These are not mutual Jacobi transformations.)

Over a period of several hours, an automated grid search {\it cum} Newton's method
turned up the following approximate $(y,q)$ pairs for $s_n(1)$:
$$ \matrix{y&q\cr\cr 1.2482046102601\,i+0.04657952537373&%
 0.2041895179564\,i+0.06533908137423\cr 0.78322226431624\,i+%
 0.30719109513916&0.45270853094805\,i+0.20573877584467\cr %
 0.49235539271999-0.74875275029651\,i&0.69018582634555\,i+%
 0.54229640598463\cr 1.0429573809123\,i+1.37265463822724&%
 0.73581195373709\,i+0.23968436737044\cr 1.00790006282379\,i+%
 0.75028009770783&0.73581195373709\,i-0.23968436737042\cr
 1.11453161008963\,i+0.62943384983216&0.63418111840451\,i-%
 0.43035475675355\cr 0.2917595098002\,i+1.21452080661459&%
 0.2041895179564\,i+0.06533908137423\cr\cr %
 0.4615909499519\,i+1.463034339711&0.86308985589011\,i+0.319843560314%
 \cr 2.01129861135631\,i+1.12816860769464&0.45270853094804\,i+%
 0.20573877584469\cr 0.23882852142275\,i+%
 1.76522302148841&0.73581195373708\,i+0.23968436737043\,.\cr} $$
The last three are quite unlike the first seven for noninteger $n$.
For integer $n$, $s_n(1)$ is real and the conjugates of all these work as well.
Note the equality of the first and last $q$ in the first group.  This would
seem to be an instance of translation by a quasiperiod, via
$$ \vartheta_{1}\left(y+i\,n\,\log q,q\right)={{\left(-1\right)^{n}\,%
 e^{2\,i\,n\,y}\,\vartheta_{1}\left(y,q\right)}\over{q^{n^{2}}}} $$
for some integer $n$.  But there are three problems with this.  First,
following a tradition that still puzzles me, we have made no provision
for a geometric ($r^n$) factor in our $\vartheta$ formul\ae, even though
the recurrence relation is unaffected by such a factor.  (Ah, but the EDS
 relation {\it is} affected.)  Second, if we
go ahead and solve for $n$,
$$ n=0.75842109957414\,i. $$
Not an integer, but pure imaginary, for some reason.  Third, if, for some
integer $n$, one of these turned up in our search, why wouldn't we find
dozens more engendered by other values of $n$?

Also compare
the fourth and fifth $q$ of the first group with the last $q$ of the second
group.  This offers hope for some simple relation between the corresponding $y$.

The grid search also turned up three spurious pairs,
$$ \matrix{y&q\cr\cr 0.16973145507896\,i+0.59439464562658&%
 0.7091459745365\,i+0.51937050102005\cr 0.15670471616132\,i+%
 0.57764641771107&0.61450698230835\,i-0.63841374520714\cr %
 0.13574561978509\,i+3.27166894673346&0.39191084430236\,i+%
 0.91218127829481,\cr } $$
which generate the sequence
\begin{eqnarray*}r_{n}&=&-r_{-n}={{144\,r_{n-3}\,r_{n-1}+432\,r^{2}_{n-2}}\over{r_{n-4}}}\cr
 &=&\ldots, -1,0,1,2^{2}\,3,-2^{4}\,3^{3},2^{7}\,3^{6},2^{12}\,3^{10},%
 -2^{16}\,3^{15},-2^{23}\,3^{20},-2^{29}\,3^{26}\,5,2^{37}\,3^{33}\,7%
 ,2^{46}\,3^{41},-2^{55}\,3^{50}\,13,\\
& &\qquad2^{68}\,3^{61},2^{77}\,3^{70}\,%
 31,2^{90}\,3^{81}\,29,-2^{103}\,3^{93}\,181,-2^{117}\,3^{106}\,5\,53%
 ,2^{133}\,3^{120}\,11\,17,-2^{148}\,3^{135}\,7\,107,\,\ldots,\\ \end{eqnarray*}
whose significance thus far eludes me, although it appears to be a (weak) EDS.


{\bf Somos5} differs from Somos4 in two respects:
$b_n=b_{4-n}$ instead of $a_n=a_{3-n}$, and a different order-reducing
relation from Conjecture 4a.

Substituting into Conjecture 4.5 $s_n=b_n, u=x, y=0, z=-1, v=1, w=0,$
\begin{eqnarray*}D_b\pmatrix{x+1/2,&3/2,&1/2\cr x-1/2,&-1/2,&-3/2\cr}&=&
\left| \matrix{b_{1}\,b_{2\,x}&b_{x}\,b_{x+1}&b_{x-1}\,b_{x+2}\cr b_{2-x%
 }\,b_{x+1}&b_{1}\,b_{2}&b_{0}\,b_{3}\cr b_{1-x}\,b_{x}&b_{0}\,b_{1}&%
 b_{-1}\,b_{2}\cr}\right|\\
\noalign{\vskip5pt}
&=&\left|\matrix{b_{2\,x}&b_{x}\,b_{x+1}&b_{x-1}\,b_{x+2}%
 \cr b_{x+1}\,b_{x+2}&1&1\cr b_{x}\,b_{x+3}&1&2\cr }\right|\;\;=\;\;0.\cr \end{eqnarray*}
This gives us $b_{2x}$ in terms of five values near $b_x$.
Alternatively, put $x-1$ for $x$, then $u=x, y=2, z=0, v=0, w=1$ to get
\begin{eqnarray*}D_b\pmatrix{x+1/2,&1/2,&3/2\cr x-3/2,&3/2,&-1/2\cr}&=&
\left|\matrix{b_{2}\,b_{2\,x-1}&b_{x-1}\,b_{x+2}&b_{x}\,b_{x+1}\cr b_{2%
 -x}\,b_{x-1}&b_{-1}\,b_{2}&b_{0}\,b_{1}\cr b_{3-x}\,b_{x}&b_{0}\,b_{%
 3}&b_{1}\,b_{2}\cr}\right|\\
\noalign{\vskip5pt}
 &=&\left|\matrix{b_{2\,x-1}&b_{x-1}\,b_{x+2}&b_{x}\,b_{%
 x+1}\cr b_{x-1}\,b_{x+2}&2&1\cr b_{x}\,b_{x+1}&1&1\cr }\right|\;\;=\;\;0\rm.\cr
\end{eqnarray*}

This gives us $b_{2x-1}$ in terms of four consecutive $b$ values.  
However, we can reduce these to {\it three}!
$$D_b\pmatrix{k-1/2,&k+1/2,&k+3/2\cr k-1/2,&k+1/2,&k-3/2\cr}=
\left| \matrix{b_{2\,k-1}&2\,b_{2\,k}&b_{2\,k-2}\cr b_{2\,k}&b_{2\,k+1}&%
 b_{2\,k-1}\cr b_{2\,k+1}&b_{2\,k+2}&b_{2\,k}\cr }\right|=0\rm, $$
or
$$ b_{2\,k+2}=-{{b_{2\,k-2}\,b^{2}_{2\,k+1}-3\,b_{2\,k}\,b_{2\,k-1}\,%
 b_{2\,k+1}+2\,b^{3}_{2\,k}}\over{b^{2}_{2\,k-1}-b_{2\,k}\,b_{2\,k-2}%
 }}\rm, $$
a fourth order recurrence.  Alternatively,
$$\twoline{D_b\pmatrix{k-3/2,&k-1/2,&k+1/2\cr k-1/2,&k+1/2,&k+3/2\cr}=}{5pt}
{ -2\,\left(\left(3\,b^{2}_{2\,k-1}-4\,b_{2\,k}\,b_{2\,k-2}\right)\,%
 b_{2\,k+2}+3\,b_{2\,k-2}\,b^{2}_{2\,k+1}-7\,b_{2\,k}\,b_{2\,k-1}\,b%
 _{2\,k+1}+5\,b^{3}_{2\,k}\right) }$$
or
$$ b_{2\,k+2}=-{{3\,b_{2\,k-2}\,b^{2}_{2\,k+1}-7\,b_{2\,k}\,b_{2\,k-1%
 }\,b_{2\,k+1}+5\,b^{3}_{2\,k}}\over{3\,b^{2}_{2\,k-1}-4\,b_{2\,k}\,b%
 _{2\,k-2}}}\rm, $$
a {\it different} fourth order recurrence.  Subtracting,
$$ b^{2}_{2\,k-2}\,b^{2}_{2\,k+1}+\left(2\,b^{3}_{2\,k-1}-5\,b_{2\,k}%
 \,b_{2\,k-2}\,b_{2\,k-1}\right)\,b_{2\,k+1}-b^{2}_{2\,k}\,b^{2}_{2\,%
 k-1}+3\,b^{3}_{2\,k}\,b_{2\,k-2}\rm, $$
a third order recurrence for $b_{2k+1}$ in terms of the three previous terms.
(Assuming you know which sign to take on the square root).  But what about
$b_{2k}$?  Simply replace $k$ by $3-k$ and $b_x$ by $b_{4-x}$, and we have
$b_{2k}$ in terms of $b_{2k-1}, b_{2k-2},$ and $b_{2k-3}$:
$$ 2\,b_{2\,k-3}\,b^{3}_{2\,k-1}-b^{2}_{2\,k-2}\,b^{2}_{2\,k-1}-5\,b%
 _{2\,k}\,b_{2\,k-3}\,b_{2\,k-2}\,b_{2\,k-1}+3\,b_{2\,k}\,b^{3}_{2\,k%
 -2}+b^{2}_{2\,k}\,b^{2}_{2\,k-3}=0 $$

But this sensitivity $\mod 2$ entails four residue classes when we order-reduce
the duplication formul\ae.:
\begin{eqnarray*}0&=&
 b^{2}_{2\,x+2}\,b^{2}_{4\,x-1}+6\,b^{3}_{2\,x}\,b^{3}_{2\,x+2}\,b%
 _{4\,x-1}-21\,b^{2}_{2\,x}\,b^{2}_{2\,x+1}\,b^{2}_{2\,x+2}\,b_{4\,x-%
 1}+16\,b_{2\,x}\,b^{4}_{2\,x+1}\,b_{2\,x+2}\,b_{4\,x-1}\\
  & &\qquad-4\,b^{6}_{2%
 \,x+1}\,b_{4\,x-1}+9\,b^{6}_{2\,x}\,b^{4}_{2\,x+2}-36\,b^{5}_{2\,x}%
 \,b^{2}_{2\,x+1}\,b^{3}_{2\,x+2}+57\,b^{4}_{2\,x}\,b^{4}_{2\,x+1}\,b%
 ^{2}_{2\,x+2}-36\,b^{3}_{2\,x}\,b^{6}_{2\,x+1}\,b_{2\,x+2}\\
   & &\qquad+8\,b^{2}%
 _{2\,x}\,b^{8}_{2\,x+1} \\
\noalign{\vskip3pt}
&=& -3\,b^{2}_{2\,x}\,b^{2}_{2\,x+2}+4\,b_{2\,x}\,b^{2}_{2\,x+1}\,b_{2%
 \,x+2}-2\,b^{4}_{2\,x+1}+b_{4\,x} \\
\noalign{\vskip3pt}
&=& b^{2}_{2\,x}\,b^{2}_{4\,x+1}+6\,b^{3}_{2\,x}\,b^{3}_{2\,x+2}\,b_{4%
 \,x+1}-21\,b^{2}_{2\,x}\,b^{2}_{2\,x+1}\,b^{2}_{2\,x+2}\,b_{4\,x+1}+%
 16\,b_{2\,x}\,b^{4}_{2\,x+1}\,b_{2\,x+2}\,b_{4\,x+1}\\
   & &\qquad-4\,b^{6}_{2\,x+%
 1}\,b_{4\,x+1}+9\,b^{4}_{2\,x}\,b^{6}_{2\,x+2}-36\,b^{3}_{2\,x}\,b^{%
 2}_{2\,x+1}\,b^{5}_{2\,x+2}+57\,b^{2}_{2\,x}\,b^{4}_{2\,x+1}\,b^{4}%
 _{2\,x+2}-36\,b_{2\,x}\,b^{6}_{2\,x+1}\,b^{3}_{2\,x+2}\\
   & &\qquad+8\,b^{8}_{2\,%
 x+1}\,b^{2}_{2\,x+2} \\
\noalign{\vskip3pt}
&=& b^{4}_{2\,x}\,b^{2}_{4\,x+2}-6\,b^{4}_{2\,x}\,b^{4}_{2\,x+2}\,b_{4%
 \,x+2}+32\,b^{3}_{2\,x}\,b^{2}_{2\,x+1}\,b^{3}_{2\,x+2}\,b_{4\,x+2}-%
 62\,b^{2}_{2\,x}\,b^{4}_{2\,x+1}\,b^{2}_{2\,x+2}\,b_{4\,x+2}\\
   & &\qquad+40\,b_{%
 2\,x}\,b^{6}_{2\,x+1}\,b_{2\,x+2}\,b_{4\,x+2}-8\,b^{8}_{2\,x+1}\,b_{%
 4\,x+2}+9\,b^{4}_{2\,x}\,b^{8}_{2\,x+2}-48\,b^{3}_{2\,x}\,b^{2}_{2\,%
 x+1}\,b^{7}_{2\,x+2}\\
   & &\qquad+86\,b^{2}_{2\,x}\,b^{4}_{2\,x+1}\,b^{6}_{2\,x+2%
 }-56\,b_{2\,x}\,b^{6}_{2\,x+1}\,b^{5}_{2\,x+2}+12\,b^{8}_{2\,x+1}\,b%
 ^{4}_{2\,x+2}\rm.\\ \end{eqnarray*}

We can obviate the first or last of these with (respectively) the odd or even
version of the third order recurrence.  This takes care of doubling.

As with Somos4, we assume three or four values near $b_{nx}$ and another tuple
near $b_{(n+1)x}$.  Then

\begin{eqnarray*}D_b\pmatrix{(n+1)x+1/2,&1/2,&3/2\cr nx-1/2,&-3/2,&-1/2}&=&\left|\matrix{b_{x+1%
 }\,b_{\left(2\,n+1\right)\,x}&b_{\left(n+1\right)\,x-1}\,b_{\left(n+%
 1\right)\,x+2}&b_{\left(n+1\right)\,x}\,b_{\left(n+1\right)\,x+1}%
 \cr b_{n\,x}\,b_{n\,x+3}&2&1\cr b_{n\,x+1}\,b_{n\,x+2}&1&1\cr }\right|\\
\noalign{\vskip3pt}
   &=&0\\ \end{eqnarray*}

gives $b_{(2n+1)x}$.  Similar constructions provide the adjacent values, and
in principle, we can use the third order relations to make everything work on
triples.

Note, however, that we could avoid the square roots and $\mod 4$ intricacies by
maintaining four values instead of three, with the help of the fourth order
relations that we subtracted to get the third order one.

(Brief flame:  a nearly forgotten fact of hardware design is that a binary
square root instruction via the ``schoolboy algorithm'' is actually simpler
than the divide instruction.  In the early 1960s, the Packard-Bell 250,
as feeble a machine as you could imagine, whose active registers were
magnetostrictive delay lines instead of flip-flop words, and whose divide
instruction needed a software followup correction, nevertheless had a hardware
square root (with remainder) that worked perfectly, in the same time as an
uncorrected divide.)

But which sign of the square root do we take?  Not obvious!  {\it E.g.},
suppose we try to use the third order Somos4 relation to compute $a_x$ from the
three previous values:
$$ a_{x}=s_x{{\sqrt{-4\,a^{3}_{x-3}\,a^{3}_{x-1}+12\,a^{2}_{x-3}\,a^{2}%
 _{x-2}\,a^{2}_{x-1}-8\,a_{x-3}\,a^{4}_{x-2}\,a_{x-1}+a^{6}_{x-2}}%
 }\over{2\,a^{2}_{x-3}}}+{{2\,a_{x-2}\,a_{x-1}}\over{a_{x-3}}}-{{a^{3%
 }_{x-2}}\over{2\,a^{2}_{x-3}}}\rm. $$

Then for $2\le x\le 38$, the sign $s_x$ coincides with
$$ {\rm sgn}\left(\phi\,\left(x-1\right)-{\rm rnd}\left(\phi\,%
 \left(x-1\right)\right)\right)\rm, $$
where $\phi$ is the golden ratio, and rnd$(x):=\lfloor x+1/2\rfloor$, the ``round''
 function.  But for $x=39$, this fails!

In practice this isn't really a problem, since we can simply choose
whichever produces an integer value for $a_x$, and we can usually check this
modulo something small.  But there may be another solution.

To take a simpler example, from the third order relation for $b_n$, we notice that
$$ \sqrt{b_{2\,n+1}\,b_{2\,n+3}-b^{2}_{2\,n+2}} 
=1, 0, 1, 1, 8, 57, 455, 22352, 47767, 69739671,\ldots$$
may be an elliptic divisibity sequence.  Trying various sign patterns, we
eventually find the recurrence
$$ h_{n}=-{{8\,h_{n-4}\,h_{n-1}+57\,h_{n-3}\,h_{n-2}}\over{h_{n-5}}} 
=1, 0,  - 1, 1, 8, 57,  - 455, 22352, 47767, 69739671, 3385862936,\ldots .$$
But this disobeys the EDS formula.  Searching further,
$$ g_{n}={{57\,g_{n-3}\,g_{n-2}-8\,g_{n-4}\,g_{n-1}}\over{g_{n-5}}}
= -1,0,1,-1,-8,57,455,-22352,-47767,69739671,-3385862936,\ldots  $$
appears to satisfy the EDS addition formula.  Evidently, $g_n=-h_n$ except when
$n=(0\mod4)$.

It seems reasonable to conjecture that the desired sign pattern for these square
root expressions is one that yields an EDS, or at least a simple recurrence.  An
example of the latter is

$$ \sqrt{3\,b_{2\,n}\,b_{2\,n+2}-2\,b^{2}_{2\,n+1}}=1,1,1,1,7%
 ,1,391,2729,175111,\ldots\rm,  $$
which cannot be an EDS.  But there is an assignation of signs:

$$ f_{n}={{57\,f_{n-3}\,f_{n-2}-8\,f_{n-4}\,f_{n-1}}\over{f_{n-5}}}=%
 -1, -1, 1, 1,-7,-1,391,-2729,-175111,8888873,565353361,\ldots\rm , $$
{\it i.e.}, the same recurrence as $g_n$, above.

Unfortunately, the order of the recurrence is likely to exceed the order of the
relation that engendered the square root, defeating the presumable purpose of
maintaining smaller intervals of consecutive values.

However, these $f$ and $g$ sequences serve another purpose if we interlace them:
\begin{eqnarray*}B_n&=&{{B_{n-4}\,B_{n-1}+B_{n-3}\,B_{n-2}}\over{B_{n-5}}}\;\;=\;\;-B_{-n}\\
&=&-1,0,1,1,1,-1,-7,-8,-1,57,391,455,-2729,-22352,-175111,-47767,8888873,\ldots,\\ \end{eqnarray*}
which is not quite an EDS.  Yet we have
$$ b_{k}={{B_{2n+1}\,b_{k-1}\,b_{k-2n-2}-B_{2n+2}\,b%
 _{k-2}\,b_{k-2n-1}}\over{B_{2n}b_{k-2n-3}}}, $$
generalizing the Somos5 recurrence.

Note that $B_n$ obeys the Somos5 recurrence, yet somehow jumps from 1 to 7, via
the well known identity $-7=0/0$.  (Of course, in the limit of absurdity, any sequence
which is alternately 0 satisfies Somos5 and Somos7.)

The $k$-tuple speedup formula is
$$ B_{k}\,B_{2\,k}\,b_{n}\,b_{n+5\,k}=B_{2\,k}\,B_{3\,k}\,b_{n+k}\,b%
 _{n+4\,k}-B_{k}\,B_{4\,k}\,b_{n+2\,k}\,b_{n+3\,k} .$$

A three-variable generalization of these last two relations is
$$ b_{n}={{B_{2\,j}\,B_{2\,k+j}\,b_{n-j}\,b_{n-2\,k-2\,j}-B_{j}\,B_{2%
 \,k+2\,j}\,b_{n-2\,j}\,b_{n-2\,k-j}}\over{B_{j}\,B_{2\,k}\,b_{n-2\,k%
 -3\,j}}}. $$

This is the same relation as with $a$ and $A$, except that $k$ must be even.
The further generalization $k \leftarrow k/2$ involves an additional term that
appears guessable, and seems to vanish for even $j$.

Thus, if $b_n$ is the Somos5 $\cosh$, then $B_n$ is the $\sinh$, although
a glitch is that $b_n$ is centered at $n=2$ rather than $n=0$.  As with $a$ and
$A$, we can merely substitute $B$ for $b$ in the last identity:

$$ B_{j}\,B_{2\,k}\,B_{n}\,B_{n-2\,k-3\,j}=B_{2\,j}\,B_{2\,k+j}\,B_{n%
 -j}\,B_{n-2\,k-2\,j}-B_{j}\,B_{2\,k+2\,j}\,B_{n-2\,j}\,B_{n-2\,k-j}. $$

This all suggests that $A$ and $B$ will have cheaper addition algorithms than
$a$ and $b$ (Somos4 and 5).  But caution:  even though $B_n$ isn't quite an EDS,
its values $\mod 19$ lie in $\{0,1,7,8,11,12,18\}$.  There may be other moduli
with even (proportionately) sparser residue classes.

As with $A_4$, we can replace the $B_5=-7\ (=0/0)$ term by $x$ to get a sequence
of polynomials:

 $$ \twoline{s_n=-1,0, 1,1,1,-1,x,x-1,-1,x^{2}-x+1,-x^{3}+x^{2}-1,-x\,\left(x^{2}%
 -2\,x+2\right),-x^{4}+x^{3}-2\,x+1,}{5pt}{\left(x-1\right)\,\left(x^{4}-x^{%
 3}+x^{2}+1\right),-x^{6}+3\,x^{5}-3\,x^{4}+3\,x^{2}-2\,x+1,2\,x^{5}-%
 5\,x^{4}+6\,x^{3}-2\,x^{2}-x+1,\ldots.}$$

These do not form an EDS, even with reassignation of signs, except when $x=-2$.  But
they retain ``strong'' divisibility,
$$(s_k(x),s_n(x))=|s_{(k,n)}(x)|,$$
for all $x$.

What {\it does} (apparently) form a (weak) EDS is the (algebraic) sequence
$$ t_{n} := s_{n}\,\left(-x-1\right)^{{{\left(-1\right)^{n}+n^{2}%
 }\over{8}}}, $$

which could be made (peculiar) polynomials with $x=-1-y^8$:
$$ \displaylines{\ldots, -1,0,1,y^{5},y^{8},-y^{17},-y^{24}\,\left(y^{8}+1\right),-%
 y^{37}\,\left(y^{8}+2\right),-y^{48},y^{65}\,\left(y^{16}+3\,y^{8}+3%
 \right),\cr
 y^{80}\,\left(y^{24}+4\,y^{16}+5\,y^{8}+1\right),y^{101}\,%
 \left(y^{8}+1\right)\,\left(y^{16}+4\,y^{8}+5\right),\cr
 -y^{120}\,\left(y^{32}+5\,y^{24}+9\,y^{16}+5\,y^{8}-1\right),\ldots .}  $$

This EDS satisfies both the Somos5oid
$$ t_{n}={{t_{n-3}\,t_{n-2}\,\left(-x-1\right)^{\mathchoice {{3%
 }\over{2}}{{3}\over{2}}{3/2}{3/2}}+t_{n-4}\,t_{n-1}\,\left(-x-1%
 \right)}\over{t_{n-5}}}, $$

and the Somos4oid
$$ t_{n}={{t_{n-3}\,t_{n-1}\,\left(-x-1\right)^{\mathchoice {{5%
 }\over{4}}{{5}\over{4}}{5/4}{5/4}}-t^{2}_{n-2}\,\left(-x-1\right)%
 }\over{t_{n-4}}}.$$

Eliminating $t_{n-5}$ and $t_{n-4}$,
$$ \displaylines{\left(t^{3}_{n-2}\,t_{n}+t_{n-3}\,t^{3}_{n-1}\right)\,\left(-x-1%
 \right)^{\mathchoice {{5}\over{4}}{{5}\over{4}}{5/4}{5/4}}+t_{n-3}\,%
 t_{n-2}\,t_{n-1}\,t_{n}\,\sqrt{-x-1}+t^{2}_{n-3}\,t^{2}_{n}\cr
=\;\;t_{n-3}%
 \,t_{n-2}\,t_{n-1}\,t_{n}\,\left(-x-1\right)^{\mathchoice {{3}\over{%
 2}}{{3}\over{2}}{3/2}{3/2}}+t^{2}_{n-2}\,t^{2}_{n-1}\,\left(-x-1%
 \right).} $$

Solving for $t_n$ yields a radical whose sign seems to depend on $x$.

In the special case $x^3+5x^2-10x+11=0$, this EDS has the elementary formula
$$t_n=\left({41-35x-x^2\over23}\right)^{(n-1)(n+1)\over8}U_{n-1}\left({1\over2}
\root8\of{11x^2+17x-244\over23}\right).$$

This is basis for the Chebychev expansion for the ``sinh'' analog of Somos5.

As in the Somos4 analog, the $s_{\rm prime}$ polynomials are irreducible, at least through
$s_{67}$.

The polynomial degrees go

$$0, 0, 0, 0, 1, 1, 0, 2, 3, 3, 4, 5, 6, 5, 8, 9, 10, 11, 13, 14, 14, 17, 19, 20, 22, 24, 26, 26,\ldots,$$

being fourteen interlaced quadratic progressions:

$$\deg s_{14q+r}=(7q+r)q+[-2,0, 0, 0, 0, 1, 1, 0, 2, 3, 3, 4, 5, 6]_r,\qquad
 0\le r\le 13.$$

The polynomials appear to be monic, except for $s_{7n}$, whose leading coefficients
appear to be $\pm n$.

Happily, the three-variable relation seems to hold for general $x$:
$$ s_{j}\,s_{2\,k}\,s_{n}\,s_{n-2\,k-3\,j}=s_{2\,j}\,s_{2\,k+j}\,s_{n%
 -j}\,s_{n-2\,k-2\,j}-s_{j}\,s_{2\,k+2\,j}\,s_{n-2\,j}\,s_{n-2\,k-j}. $$

Better yet, putting $u=-x, v=-y$ and $s_n=-s_{-n}$ in Conjecture 4.5, we get
the four (integer) variable, three term relation
$$ s_{k-j+1}\,s_{k+j}\,s_{-n-i+1}\,s_{n-i}=s_{-j-i+1}\,s_{j-i}\,s_{-n%
 +k+1}\,s_{n+k}+s_{k-i+1}\,s_{k+i}\,s_{-n-j+1}\,s_{n-j}. $$

Caution:  This fails for noninteger $i,k,j,n$ even though the subscripts are
integral.

In his email to sci.math, Elkies makes the remarkable observation (modulo typos)
that $t_n:=(2/3)^{(n \mod 2)/4}b_n$ satisfies the reduced order (quasiSomos4) recurrence
$$t_{n-2}t_{n+2}=\sqrt6\,t_{n-1}t_{n+1}-t_n^2\rm\eqno{\rm(Elkies)}.$$
It is probable that $t_n$ also satisfies a third order relation (of higher degree).

As with the Somos4 polynomials, it appears that $s_n(\alpha)$ falls into $k$ interlaced
elementary progressions when $s_k(\alpha)=0$, but they are more complicated.  Alternatively,
they can be written as $mk$ simpler progressions, for some multiple $m$.  {\it E.g.},
$s_8(e^{i\pi/3})=0$ and $s_n(e^{i\pi/3})$ is merely periodic, but the period is forty-eight!

Also like the Somos4 polynomials, we can get a $\vartheta_1$ expression via the change of
variables:

$$ t_{n}(x)={{\tan \left(\arctan \root 8 \of{-x-1}+{{\pi\,n%
 }\over{2}}\right)}\over{\root 8 \of{-x-1}}}s_{n}(x), $$
where $t_n(x)$ satisfies \
$$ t_{n-4}\,t_{n}=t_{n-3}\,t_{n-1}\,\sqrt{-x-1}-t^{2}_{n-2}. $$

It shouldn't be hard to find a Chebychev formula for some algebraic $x$, and hence
elementary expansions for the $\vartheta$ parameters, as we did for Somos4.

{\bf Expression as $\vartheta$:}  Elkies' email gives
$$b_n=(3/2)^{(n \mod 2)/4}b\,u^{\left(n-2\right)^{2}}\,\sum_{k=-\infty }^{\infty }{q^{k^{2}}\,z^{k\,\left(n%
 -2\right)}} \rm,$$
with $q = 0.02208942811097933557356088\ldots,$
$z = 0.1141942041600238048921321\ldots,$ \\
$b = 0.9576898995913810138013844\ldots,$ and $u = 0.7889128685374661530379575\ldots.$
These constants do not appear to be in Plouffe's collection.

Similarly,
$$ s_{n}=B_{n}\,\tan \left({{\pi\,n}\over{2}}+\arctan \root 8 \of{6} \right) $$
also satisfies the (Elkies) recurrence, giving

\begin{eqnarray*} B_{n}&=&b\,\cot \left({{\pi\,n}\over{2}}+\arctan \root 8 \of{6}%
 \right)\,u^{n^{2}}\,\sum_{k=0}^{\infty }{\left(-1\right)^{k}\,q^{%
 \left(2\,k+1\right)^{2}}\,\sin \left(\left(2\,k+1\right)\,n\,y%
 \right)}\cr
 &=&{{b}\over{2}}\,\cot \left({{\pi\,n%
 }\over{2}}+\arctan \root 8 \of{6}\right)\,u^{n^{2}}\,\vartheta_{1}%
 \left(n\,y,q^{4}\right),\cr\end{eqnarray*}

where
$$ b=-1.82905778669392,\ u=-1.07425451486466,\ y=0.89396990235568,\ %
 q=-0.52353014451686. $$

Numerically testing this equation for $-1\le n\le 22$:

$ - 1 =  - 0.99999999999986d0,\  0 = 0.0d0,\  1 = 0.99999999999986d0, \  1 = 0.99999999999986d0, \\ \  
1 = 0.99999999999984d0,\   - 1 =  - 0.99999999999984d0,  - 7 =  - 6.99999999999914d0,
 - 8 =  - 7.99999999999886d0,\\ \   - 1 =  - 0.99999999999982d0,\  57 = 56.9999999999918d0,\  
391 = 390.999999999952d0,\  455 = 454.999999999928d0, \\ \   - 2729 =  - 2728.99999999958d0, \  
 - 22352 =  - 22351.9999999972d0,\   - 175111 =  - 175110.999999977d0, \\ \  
 - 47767 =  - 47766.9999999918d0,\  8888873 = 8888872.99999885d0,\  
69739671 = 6.9739670999992d+7, \\ \  565353361 = 5.65353360999916d+8,\  
 - 3385862936 =  - 3.3858629359995d+9, \\ \   - 195345149609 =  - 1.9534514960898d+11,\  
 - 1747973613295 =  - 1.74797361329478d+12, \\ \   - 4686154246801 =  - 4.6861542468004d+12,\  
632038062613231 = 6.32038062613152d+14$.

In email to sci.math, Randall Rathbun and Ralph Buchholz make the remarkable claim
that the Heron triangles with two rational medians have side lengths
$$\displaylines{\left[ B_{i+3}\,\left(B^{2}_{i}\,b^{2}_{i+3}\,b^{4}_{i+4}+B^{2}_{%
 i+1}\,b^{2}_{i+2}\,B^{4}_{i+2}\right)\,b_{i+5},\right.\cr
 B_{i+2}\,b_{i+4}\,%
 \left(B^{2}_{i}\,b^{2}_{i+2}\,B^{2}_{i+3}\,b^{2}_{i+5}+B^{2}_{i+1}\,%
 B^{2}_{i+2}\,b^{2}_{i+3}\,b^{2}_{i+4}\right),\cr
 \left.B_{i+1}\,b_{i+3}\,%
 \left(B^{2}_{i}\,4^{i+1\mod2}\,B^{4}_{i+2}\,b%
 ^{2}_{i+5}+4^{i\mod2}\,b^{2}_{i+2}\,B^{2}_{i+3%
 }\,b^{4}_{i+4}\right)\right]. \cr} $$


{\bf Somos6:}
\begin{eqnarray*} D_c\pmatrix{n-3,&0,&1,&2\cr 0,&1,&2,&3\cr }&=&
\left|\matrix{c^{2}_{n-3}&c_{n-4}\,c_{n-2}&c_{n-5}\,c_{n-1}&c_{n-6}\,c%
 _{n}\cr 1&3&5&9\cr 1&1&3&5\cr 1&1&1&3\cr }\right|\\
\noalign{\vskip3pt}
&=&-4\,c_{n-6}\,c_{n}+4\,c_{ n-5}\,c_{n-1}+4\,c_{n-4}\,c_{n-2}+4\,c^{2}_{n-3},\\ \end{eqnarray*}

(four times) the defining recurrence for Somos6.  But
 
$$D_c \pmatrix{0,&2,&4,&6\cr 0,&1,&3,&4\cr }=80=D_c\pmatrix{-\mathchoice {{1}\over{2}}{{1}\over{2}}{1/2}{1/2},&%
 \mathchoice {{1}\over{2}}{{1}\over{2}}{1/2}{1/2},&\mathchoice {{3%
 }\over{2}}{{3}\over{2}}{3/2}{3/2},&\mathchoice {{13}\over{2}}{{13%
 }\over{2}}{13/2}{13/2}\cr \mathchoice {{1}\over{2}}{{1}\over{2}}{1/2%
 }{1/2},&\mathchoice {{3}\over{2}}{{3}\over{2}}{3/2}{3/2},&\mathchoice %
 {{7}\over{2}}{{7}\over{2}}{7/2}{7/2},&\mathchoice {{5}\over{2}}{{5%
 }\over{2}}{5/2}{5/2}\cr }, $$
so 4 by 4 isn't enough.  Building on the nonsingular matrix,
$$D_c\pmatrix{x,&0,&2,&4,&6\cr y,&0,&1,&3,&4\cr }=
\left|\matrix{c_{x-y}\,c_{y+x}&c^{2}_{x}&c_{x-1}\,c_{x+1}&c_{x-3}\,c_{x%
 +3}&c_{x-4}\,c_{x+4}\cr c_{-y}\,c_{y}&1&3&9&23\cr c_{2-y}\,c_{y+2}&1%
 &1&3&15\cr c_{4-y}\,c_{y+4}&1&1&5&9\cr c_{6-y}\,c_{y+6}&9&5&23&75%
 \cr }\right|$$
$$\displaylines{ =80\,c_{x-y}\,c_{y+x}+\left(4
 \,c_{x-4}\,c_{x+4}+12\,c_{x-3}\,c_{x+3}-52\,c_{x-1}\,c_{x+1}-44\,c^{%
 2}_{x}\right)\,c_{y-1}\,c_{y+6}\cr
+\left(-16\,c_{x-4}\,c_{x+4}-88\,c_{x%
 -3}\,c_{x+3}+328\,c_{x-1}\,c_{x+1}+176\,c^{2}_{x}\right)\,c_{y+1}\,c%
 _{y+4}\cr
+\left(-28\,c_{x-4}\,c_{x+4}-44\,c_{x-3}\,c_{x+3}+284\,c_{x-1}%
 \,c_{x+1}+188\,c^{2}_{x}\right)\,c_{y+3}\,c_{y+2}\cr
+\left(8\,c_{x-4}\,c_{x+4}+24\,c_{x-3}\,c_{x+3}-144\,c_{x-1}\,c_{x+1}-48\,c^{2}_{x}%
 \right)\,c_{y+5}\,c_{y},}$$
which empirically vanishes, providing a fairly messy addition formula.  Using the
defining recurrence to algebraically eliminate $c_{y+6}, c_{x-4},$ etc., yields
even messier relations, but with lower order ($=$ width).  The determinant approach
can't do much better unless we can find a ``narrower'' non-singular 4 by 4.  But this is
unlikely, since replacing the sixteen numerical coefficients by undetermined ones
indicates that the $c_{y-1}c_{y+6}$ term is indispensable, and the $c_{x-4}c_{x+4}$
term is not replaceable by $c_{x-2}c_{x+2}$.  Slightly better may be

$$ \displaylines{80\,c_{x-y+1}\,c_{y+x}=\left(-4\,c_{x-3}\,c_{x+4}-8\,c_{x-2}\,c_{%
 x+3}+16\,c_{x-1}\,c_{x+2}+28\,c_{x}\,c_{x+1}\right)\,c_{y-2}\,c_{y+6%
 }\cr
+\left(32\,c_{x-3}\,c_{x+4}+24\,c_{x-2}\,c_{x+3}-88\,c_{x-1}\,c_{x+%
 2}-144\,c_{x}\,c_{x+1}\right)\,c_{y-1}\,c_{y+5}\cr
+\left(-44\,c_{x-3}\,%
 c_{x+4}-48\,c_{x-2}\,c_{x+3}+176\,c_{x-1}\,c_{x+2}+188\,c_{x}\,c_{x+%
 1}\right)\,c_{y}\,c_{y+4}\cr+\left(8\,c_{x-3}\,c_{x+4}+96\,c_{x-2}\,c_{%
 x+3}-152\,c_{x-1}\,c_{x+2}-96\,c_{x}\,c_{x+1}\right)\,c_{y+1}\,c_{y+%
 3}\cr} $$
from
$$D_c\pmatrix{x+{1\over2},&-\mathchoice {{1}\over{2}}{{1}\over{2}}{1/2}{1/2},&%
 \mathchoice {{1}\over{2}}{{1}\over{2}}{1/2}{1/2},&\mathchoice {{3%
 }\over{2}}{{3}\over{2}}{3/2}{3/2},&\mathchoice {{13}\over{2}}{{13%
 }\over{2}}{13/2}{13/2}\cr
y-{1\over2},& \mathchoice {{1}\over{2}}{{1}\over{2}}{1/2%
 }{1/2},&\mathchoice {{3}\over{2}}{{3}\over{2}}{3/2}{3/2},&\mathchoice %
 {{7}\over{2}}{{7}\over{2}}{7/2}{7/2},&\mathchoice {{5}\over{2}}{{5%
 }\over{2}}{5/2}{5/2}\cr }=0. $$

Conjecture:  there is no sequence of bivariate polynomials obeying Somos6.  Evidence:
Initializing with

$$s_0,s_1,\ldots=0,1,x^{3}\,y,-x^{6}\,y^{2},x^{5}\,y^{2},-x^{8}\,y^{3}\,%
 \left(x^{4}\,y+1\right),-x^{15}\,y^{5},\ldots,$$

gives polynomials through $s_{29}$, but then $s_{30}$ has a denominator of $x$.  Replacing
$y\leftarrow yx$ will probably move the violation a few terms to the right.

However, if we initialize with
$$s_0,s_1,\ldots=0,1,1,1,1,-2,x,\ldots,$$
then we appear to get polynomials 
$$\displaylines{s_7,s_8,\ldots= x-1,2\,x+3,x^{2}+5,x^{2}+x-9,-\left(x^{3}+2\,x^{2}+4\,x+2%
 \right),x^{2}+13\,x-13,\cr
-\left(x^{4}+2\,x^{3}+7\,x^{2}+7\,x+32\right)%
 ,x^{4}-x^{3}+10\,x^{2}-41\,x-13,x^{5}+5\,x^{4}+10\,x^{3}+7\,x^{2}+32%
 \,x+70,\ldots\cr}$$
which is decidedly not a divisibility sequence.  In fact, the only reducible $s_n$ through
$s_{69}$ is 
$$ s_{16}=-\left(x+1\right)\,\left(3\,x^{4}-x^{3}+21\,x^{2}-52\,x+117\right)\qquad! $$
And, unlike Somos4 and 5, $s_k(\alpha)=0 \not \Rightarrow s_{mk}(\alpha)=0$.

{\bf Somos7:}  Random clue:
$$ d_{k-11}\,d_{k}+d_{k-10}\,d_{k-1}+d_{k-7}\,d_{k-4}=61\,d_{k-6}\,d_{k-5}. $$

\begin{eqnarray*}D_{d}\pmatrix{n-\mathchoice {{7}\over{2}}{{7}\over{2}}{7/2}{%
 7/2},&\mathchoice {{1}\over{2}}{{1}\over{2}}{1/2}{1/2},&\mathchoice {{%
 3}\over{2}}{{3}\over{2}}{3/2}{3/2},&\mathchoice {{5}\over{2}}{{5%
 }\over{2}}{5/2}{5/2}\cr \mathchoice {{1}\over{2}}{{1}\over{2}}{1/2}{%
 1/2},&\mathchoice {{3}\over{2}}{{3}\over{2}}{3/2}{3/2},&\mathchoice {{%
 5}\over{2}}{{5}\over{2}}{5/2}{5/2},&\mathchoice {{7}\over{2}}{{7%
 }\over{2}}{7/2}{7/2}\cr } &=&
\left|\matrix{d_{n-4}\,d_{n-3}&d_{n-5}\,d_{n-2}&d_{n-6}\,d_{n-1}&d_{n-7%
 }\,d_{n}\cr 1&3&5&9\cr 1&1&3&5\cr 1&1&1&3\cr }\right|\\
\noalign{\vskip3pt}
&=&-4\,\left(d_{n-7}\,d_{n}-d_{n-6}\,d_{n-1}-d_{n-5}\,d_{n-2}-d_{n-4}\,d_{n-3}\right),\cr \end{eqnarray*}
the defining recurrence.  But
$$ D_{d} \pmatrix{-\mathchoice {{3}\over{2}}{{3}\over{2}}{3/2}{3%
 /2},&\mathchoice {{1}\over{2}}{{1}\over{2}}{1/2}{1/2},&\mathchoice {{7%
 }\over{2}}{{7}\over{2}}{7/2}{7/2},&\mathchoice {{9}\over{2}}{{9%
 }\over{2}}{9/2}{9/2}\cr -\mathchoice {{3}\over{2}}{{3}\over{2}}{3/2}%
 {3/2},&\mathchoice {{1}\over{2}}{{1}\over{2}}{1/2}{1/2},&\mathchoice %
 {{5}\over{2}}{{5}\over{2}}{5/2}{5/2},&\mathchoice {{9}\over{2}}{{9%
 }\over{2}}{9/2}{9/2}\cr }=160 ,$$
so
$$ D_{d}\pmatrix{x+\mathchoice {{1}\over{2}}{{1}\over{2}}{1/2}{%
 1/2},&-\mathchoice {{3}\over{2}}{{3}\over{2}}{3/2}{3/2},&\mathchoice %
 {{1}\over{2}}{{1}\over{2}}{1/2}{1/2},&\mathchoice {{7}\over{2}}{{7%
 }\over{2}}{7/2}{7/2},&\mathchoice {{9}\over{2}}{{9}\over{2}}{9/2}{9/2%
 }\cr y-\mathchoice {{1}\over{2}}{{1}\over{2}}{1/2}{1/2},&-%
 \mathchoice {{3}\over{2}}{{3}\over{2}}{3/2}{3/2},&\mathchoice {{1%
 }\over{2}}{{1}\over{2}}{1/2}{1/2},&\mathchoice {{5}\over{2}}{{5%
 }\over{2}}{5/2}{5/2},&\mathchoice {{9}\over{2}}{{9}\over{2}}{9/2}{9/2%
 }\cr }=\left|\matrix{d_{-y+x+1}\,d_{y+x}&d_{x-1}\,d_{x+2}&d_{x}\,d_{x+1}&%
 d_{x-2}\,d_{x+3}&d_{x-4}\,d_{x+5}\cr d_{-y-1}\,d_{y-2}&9&15&17&137%
 \cr d_{1-y}\,d_{y}&3&1&5&17\cr d_{4-y}\,d_{y+3}&1&1&1&15\cr d_{5-y}%
 \,d_{y+4}&1&1&3&9\cr }\right|, $$
and
$$ \displaylines{160\,d_{x-y+1}\,d_{y+x}+\left(4\,d_{x-4}\,d_{x+5}+12\,d_{x-2}\,d%
 _{x+3}-28\,d_{x-1}\,d_{x+2}-44\,d_{x}\,d_{x+1}\right)\,d_{y-2}\,d_{y%
 +7}\cr
+\left(12\,d_{x-4}\,d_{x+5}+36\,d_{x-2}\,d_{x+3}-164\,d_{x-1}\,d%
 _{x+2}-52\,d_{x}\,d_{x+1}\right)\,d_{y}\,d_{y+5}\cr
+\left(-28\,d_{x-4}%
 \,d_{x+5}-164\,d_{x-2}\,d_{x+3}+356\,d_{x-1}\,d_{x+2}+228\,d_{x}\,d%
 _{x+1}\right)\,d_{y+1}\,d_{y+4}\cr
+\left(-44\,d_{x-4}\,d_{x+5}-52\,d_{x%
 -2}\,d_{x+3}+228\,d_{x-1}\,d_{x+2}+324\,d_{x}\,d_{x+1}\right)\,d_{y+%
 2}\,d_{y+3}\cr} $$
appears to vanish.

Messy doubling formul\ae, in case I don't find anything better in the next few days:

$$\displaylines{ d_{2\,k-1}=-(d_{k-5}\,d_{k-2}\,d_{k+4}\,d_{k+7} +3\,d_{k-3}\,d_{k-%
 2}\,d_{k+2}\,d_{k+7} -7\,d^{2}_{k-2}\,d_{k+1}\,d_{k+7}\cr
 -11\,d_{k-2}\,d%
 _{k-1}\,d_{k}\,d_{k+7} +3\,d_{k-5}\,d_{k}\,d_{k+4}\,d_{k+5} +9\,d_{k-3%
 }\,d_{k}\,d_{k+2}\,d_{k+5} -41\,d_{k-2}\,d_{k}\,d_{k+1}\,d_{k+5} -13\,%
 d_{k-1}\,d^{2}_{k}\,d_{k+5}\cr
 -7\,d_{k-5}\,d_{k+1}\,d^{2}_{k+4} -11\,d_{%
 k-5}\,d_{k+2}\,d_{k+3}\,d_{k+4} -41\,d_{k-3}\,d_{k+1}\,d_{k+2}\,d_{k+%
 4} +89\,d_{k-2}\,d^{2}_{k+1}\,d_{k+4} +57\,d_{k-1}\,d_{k}\,d_{k+1}\,d%
 _{k+4}\cr
 -13\,d_{k-3}\,d^{2}_{k+2}\,d_{k+3} +57\,d_{k-2}\,d_{k+1}\,d_{k+%
 2}\,d_{k+3} +81\,d_{k-1}\,d_{k}\,d_{k+2}\,d_{k+3})/{40}\cr} $$

$$ \displaylines{d_{2\,k}=(d_{k-6}\,d_{k-1}\,d_{k+5}\,d_{k+8} -3\,d_{k-4}\,d_{k-1}%
 \,d_{k+3}\,d_{k+8} +11\,d_{k-2}\,d_{k-1}\,d_{k+1}\,d_{k+8}\cr
 -47\,d^{2}%
 _{k-1}\,d_{k}\,d_{k+8} +3\,d_{k-6}\,d_{k+1}\,d_{k+5}\,d_{k+6} -9\,d_{k%
 -4}\,d_{k+1}\,d_{k+3}\,d_{k+6} +53\,d_{k-2}\,d^{2}_{k+1}\,d_{k+6}\cr
 -161%
 \,d_{k-1}\,d_{k}\,d_{k+1}\,d_{k+6} -7\,d_{k-6}\,d_{k+2}\,d^{2}_{k+5} -%
 11\,d_{k-6}\,d_{k+3}\,d_{k+4}\,d_{k+5} +41\,d_{k-4}\,d_{k+2}\,d_{k+3}%
 \,d_{k+5}\cr
 -137\,d_{k-2}\,d_{k+1}\,d_{k+2}\,d_{k+5} +329\,d_{k-1}\,d_{k%
 }\,d_{k+2}\,d_{k+5} +13\,d_{k-4}\,d^{2}_{k+3}\,d_{k+4} -81\,d_{k-2}\,d%
 _{k+1}\,d_{k+3}\,d_{k+4}\cr
   +577\,d_{k-1}\,d_{k}\,d_{k+3}\,d_{k+4}%
 )/{120}\cr} $$

$$ \displaylines{d_{2\,k}=(d_{k-6}\,d_{k-1}\,d_{k+5}\,d_{k+8} -3\,d_{k-3}\,d_{k-1}%
 \,d_{k+2}\,d_{k+8} +8\,d_{k-2}\,d_{k-1}\,d_{k+1}\,d_{k+8}\cr
 -50\,d^{2}_{%
 k-1}\,d_{k}\,d_{k+8} +3\,d_{k-6}\,d_{k+1}\,d_{k+5}\,d_{k+6} -9\,d_{k-3%
 }\,d_{k+1}\,d_{k+2}\,d_{k+6} +44\,d_{k-2}\,d^{2}_{k+1}\,d_{k+6}\cr
 -170\,%
 d_{k-1}\,d_{k}\,d_{k+1}\,d_{k+6} -7\,d_{k-6}\,d_{k+2}\,d^{2}_{k+5} -11%
 \,d_{k-6}\,d_{k+3}\,d_{k+4}\,d_{k+5} +41\,d_{k-3}\,d^{2}_{k+2}\,d_{k+%
 5}\cr
 -96\,d_{k-2}\,d_{k+1}\,d_{k+2}\,d_{k+5} +370\,d_{k-1}\,d_{k}\,d_{k+%
 2}\,d_{k+5} +13\,d_{k-3}\,d_{k+2}\,d_{k+3}\,d_{k+4} -68\,d_{k-2}\,d_{k%
 +1}\,d_{k+3}\,d_{k+4}\cr
   +590\,d_{k-1}\,d_{k}\,d_{k+3}\,d_{k+4})/{120}\cr} $$

As an existence proof, here is a $6^{\rm th}$ order relation:

$$\displaylines{0=\left(\left(d^{2}_{k-6}\,d_{k-5}\,d_{k-3} +d_{k-6}\,d^{2}_{k-5}\,%
 d_{k-4}\right)\,d_{k-1} +\left(d_{k-6}\,d^{2}_{k-5}\,d_{k-3} +d^{3}_{k%
 -5}\,d_{k-4}\right)\,d_{k-2} +d^{2}_{k-6}\,d^{3}_{k-3}\right.\cr
 \left.+2\,d_{k-6}\,d%
 _{k-5}\,d_{k-4}\,d^{2}_{k-3} +d^{2}_{k-5}\,d^{2}_{k-4}\,d_{k-3}%
 \right)\,d^{2}_{k} +\left(\left(d^{2}_{k-6}\,d_{k-5}\,d_{k-2}\right.\right.\cr
 \left. +d^{2}_{%
 k-6}\,d_{k-4}\,d_{k-3} +2\,d_{k-6}\,d_{k-5}\,d^{2}_{k-4}\right)\,d^{2%
 }_{k-1} +\left(2\,d_{k-6}\,d^{2}_{k-5}\,d^{2}_{k-2} +\left(2\,d^{2}_{k%
 -6}\,d^{2}_{k-3}\right.\right.\cr
 \left.-56\,d_{k-6}\,d_{k-5}\,d_{k-4}\,d_{k-3} +2\,d^{2}_{k-%
 5}\,d^{2}_{k-4}\right)\,d_{k-2} +2\,d_{k-6}\,d^{2}_{k-4}\,d^{2}_{k-3}\cr
 \left.+2\,d_{k-5}\,d^{3}_{k-4}\,d_{k-3}\right)\,d_{k-1} +d^{3}_{k-5}\,d^{3}
 _{k-2} +\left(2\,d_{k-6}\,d_{k-5}\,d^{2}_{k-3} +3\,d^{2}_{k-5}\,d_{k-4%
 }\,d_{k-3}\right)\,d^{2}_{k-2}\cr
 \left.+\left(2\,d_{k-6}\,d_{k-4}\,d^{3}_{k-3%
 } +2\,d_{k-5}\,d^{2}_{k-4}\,d^{2}_{k-3}\right)\,d_{k-2}\right)\,d_{k}%
 +\left(d^{2}_{k-6}\,d_{k-4}\,d_{k-2} +d_{k-6}\,d^{3}_{k-4}\right)\,d%
 ^{3}_{k-1}\cr
 +\left(\left(d^{2}_{k-6}\,d_{k-3} +2\,d_{k-6}\,d_{k-5}\,d_{%
 k-4}\right)\,d^{2}_{k-2} +\left(3\,d_{k-6}\,d^{2}_{k-4}\,d_{k-3} +d_{k%
 -5}\,d^{3}_{k-4}\right)\,d_{k-2}\right.\cr
 \left.+d^{4}_{k-4}\,d_{k-3}\right)\,d^{2}%
 _{k-1} +\left(\left(2\,d_{k-6}\,d_{k-5}\,d_{k-3} +d^{2}_{k-5}\,d_{k-4}%
 \right)\,d^{3}_{k-2} +\left(2\,d_{k-6}\,d_{k-4}\,d^{2}_{k-3}\right.\right.\cr
 \left.\left.+3\,d_{k-%
 5}\,d^{2}_{k-4}\,d_{k-3}\right)\,d^{2}_{k-2} +2\,d^{3}_{k-4}\,d^{2}_{%
 k-3}\,d_{k-2}\right)\,d_{k-1} +d^{2}_{k-5}\,d_{k-3}\,d^{4}_{k-2}\cr
 +2\,d%
 _{k-5}\,d_{k-4}\,d^{2}_{k-3}\,d^{3}_{k-2} +d^{2}_{k-4}\,d^{3}_{k-3}\,%
 d^{2}_{k-2}.\cr} $$

It appears that
$$s_{-1}(x),s_0(x),\ldots,\;\;=\;\; - 1, 0, 1, 1, 1, 1, 1,  - 2, x, x - 1, 2x - 3, x - 4, x^2 - 4x + 2, x^2 - x - 1, \ldots,$$
gives a (non-divisibility) sequence of polynomials.  The special case $x=1$ appears to give
eight interlaced arithmetic(!) sequences:
$$s_{8q+r}(1)\;\;=\;\;(-)^q[0,1,2q+1,1,1,1,-2q-2,1]_{0\le r\le7}.$$

A curious initialization is
$$ \displaylines{\ldots, -{{3}\over{\sqrt{2}}},1,1,1,1,1,1,-\sqrt{2},2-\sqrt{2},3-2%
 \,\sqrt{2},5-4\,\sqrt{2},10-8\,\sqrt{2},28-20\,\sqrt{2},107-76\,%
 \sqrt{2},455-322\,\sqrt{2},\ldots,}  $$
where the first term with denominator $>1$ is
$$ d_{34}={{510156039514521981558192050-360734795003990787362927953\,%
 \sqrt{2}}\over{2}}, $$
and the first term (after $d_{-1}=-3/\sqrt2$) with magnitude $>2$ is

$$d_{39}\approx 2.3813134529 .$$

{\bf References:}

Plouffe's Inverter:  http://www.lacim.uqam.ca/pi/indexf.html

Sloane's Sequence Server:  http://www.research.att.com/\Tilde njas/sequences/

Theta Functions in Macsyma \cite{RWGTHETA}

\parskip=0pt
\medskip
Zagier notes on Somos5 and elliptic curve:\par
\ \ \ http://www-groups.dcs.st-andrews.ac.uk/\Tilde john/Zagier/Problems.html\par
and \cite{Elkies}.

\medskip
$\vartheta$like double sum involving seven empirical constants for Somos6: \par
\ \ \ http://grail.cba.csuohio.edu/\Tilde somos/somos6.html .

\medskip
Jim Propp's Somos sequence site:  http://www.math.wisc.edu/\Tilde propp/somos.html .

\medskip
General:  Google search Somos sequence* and Elliptic divisibility.


\section{Modular Theta Functions}

This section explores mod $P$ analogs of some classical special functions.
We're interested in the general problem of developing modular analogs for
classical special functions of analysis.  The modular versions of exponentiation
and logarithms have been known for two centuries.  These are easily
generalized to modular trigonometric and hyperbolic functions, and
thence to elliptic functions such as sn, cn, dn, and Weierstrass $\wp$.
The modular Gudermannian function gd converts between sin and tanh.
A modular version of the amplitude function am converts between sin and sn.
I've already mentioned modular polylogarithms.
This note introduces modular theta functions.

\subsection{Elliptic Function Review}

The classical elliptic functions arose from trying to
determine the arc-length of an ellipse, by integrating
expressions involving the square roots of cubic or quartic polynomials.
Elliptic functions are analytic, complex valued functions, that take
complex number arguments.  They have two periods. One period is usually
taken to be a real number, and the other is necessarily complex.
Elliptic functions have poles as their only singularities (at finite locations).  The
two periods make a period parallelogram, covering the complex plane
in a regular pattern.  The function values repeat in each parallelogram.
There are two popular ways of discussing elliptic functions, depending
on whether the basic integral is the square root of a cubic or quartic
polynomial.  The two ways are equivalent, but the choices lead to differences
in details of the formulas.  The Jacobian elliptic functions (\cite{AbSteg,AbStegWeb}, chapter 16) lead more
naturally to theta functions.
The basic Jacobian elliptic functions are sn$(u,k)$, cn$(u,k)$, and dn$(u,k)$.
Typically $k$ is fixed in an application (it is related to the shape of
the period parallelogram), and it is often elided to simplify formulas.
These functions have two poles and two zeros in each copy of the
period parallelogram.
Some of the fundamental formulas are
\begin{eqnarray*} {\rm sn}^2u+{\rm cn}^2u & = & 1  \\
 k^2 \, {\rm sn}^2u+{\rm dn}^2u & = & 1 \\
 {\rm sn}'u & = & {\rm cn}\, u \: {\rm dn}\, u \\
 {\rm sn}(u+v) & = &
 \frac{{\rm sn}\, u \: {\rm cn} \, v \: {\rm dn} \, v + {\rm sn} \, v \: {\rm cn} \, u \: {\rm dn} \, u}
 {1 - k^2 \, {\rm sn}^2u \: {\rm sn}^2v} \\
{\rm sn} \, 2u & = & \frac{2 \, {\rm sn}\, u \: {\rm cn}\, u \: {\rm dn} \, u}{1- k^2 \, {\rm sn}^4u} \\
{\rm sn} \, u & = & {\rm sin}({\rm am} \, u)\\
\end{eqnarray*}
The sn and cn functions can be regarded as sine and cosine
of a distorted input.
The am function captures the distortion.
Chapter 16 of Abramowitz \& Stegun \cite{AbSteg,AbStegWeb} has much more,
and the classic Whittaker \& Watson \cite{ww} explains what's going on.

\subsection{Elliptic Functions vs Elliptic Curves}

An elliptic curve can be parameterized by elliptic functions,
just as a circle can be parameterized by circular (trigonometric) functions (cos and sin)
and a hyperbola by hyperbolic functions (cosh and sinh).

\begin{eqnarray*}
 Y^2 & = & k^2X^4 - (1+k^2)X^2 + 1 \\
 (X,Y) & = & ({\rm sn} \, u, \: {\rm cn} \, u \: {\rm dn} \, u) \\
 ({\rm cn} \, u \: {\rm dn} \, u)^2 & = & (1- {\rm sn}^2u)(1-k^2 {\rm sn}^2u) \\
\end{eqnarray*}

The original applications of these ideas were for real and complex applications,
but there has been a minor component of number theory ever since Diophantus (c. 200 AD) introduced
problems that reduced to cubic curves.  The theory of finding integer and
rational points on elliptic curves has undergone a major development in the
last century.  The application of elliptic curves to cryptography depends on the
fact that what works for the fields {\bf R}, {\bf C}, and {\bf Q} can often be made to work for
the fields mod $p$ and the Galois Fields GF[$p^k$].

\subsection{What's a Theta Function?}


Theta functions were probably first introduced by Euler.
They arise from some infinite products related to the partition function.
Theta functions have rapidly convergent series, and they are
closely related to elliptic functions, which makes them useful
in computing elliptic function values.

Like elliptic functions, theta functions are complex valued functions
with one complex argument, and a second shape parameter.  They only have
one true period, but they have a second quasiperiod.  The period is often taken
to be 1 or $2\pi$, while the quasiperiod is a complex number.  Together the
two define a parallelogram, as with elliptic functions.
  Changing the
argument by one period leaves the value of the theta function unchanged, while
changing the argument by the quasiperiod multiplies the theta function by
a fixed value.
In one sense, theta functions are easier than elliptic functions, since they
have only one pole and one zero per parallelogram.
There are 4 basic theta functions, with the pole located either in the corner
of the parallelogram, in the middle of a side, or in the center.
Elliptic functions are ratios of theta functions, with the same
parallelogram.

\subsection{Some Basic Properties of Theta Functions}

\begin{eqnarray*}
 \vartheta_4(z,q) & = & \sum_{n=-\infty}^{\infty} (-)^n q^{n^2} e^{2inz} \\
 \vartheta_4(z+\pi) & = & \vartheta_4(z) \\
 \vartheta_4(z+\tau) & = & e^{2i \tau} \vartheta_4(z) \\
 \vartheta_{1,2,3} & = & \vartheta_4 + {\rm half \: periods} \: \;  \pi/2, \; \tau/2, \: (\pi+\tau)/2 \\
 \vartheta_4(z,q) & = & \prod_{n=1}^{\infty} (1-q^{2n})(1 - q^{2n-1} cos \,2z + q^{4n-2}) \\
\end{eqnarray*}
$$ {\rm sn} \sim \frac{\vartheta_1}{\vartheta_4}, \;
 {\rm cn} \sim \frac{\vartheta_2}{\vartheta_4}, \;
 {\rm dn} \sim \frac{\vartheta_3}{\vartheta_4} $$
Similar products exist for $\vartheta_{1,2,3}$.
The partition generating function is $ 1/\prod_{n=1}^{\infty}(1-q^n)$.
Gosper has developed a package for manipulating theta functions in
the symbolic algebra program Macsyma \cite{RWGTHETA}.

\subsection{Theta Function Identities}

As with the elliptic functions, there is a multitude, nay, a plethora, of theta function identities.
A small sampling is presented below.
$\vartheta_i(0)$ is abbreviated to $\vartheta_i$.
\begin{eqnarray*}
 \vartheta_2^4 + \vartheta_4^4 & = & \vartheta_3 ^4 \\
 \vartheta_1^2(z) \vartheta_2^2  & = & \vartheta_3^2 \vartheta_4(z) - \vartheta_4^2 \vartheta_3(z) \\
 \vartheta_1^4(z) + \vartheta_3^4(z) & = & \vartheta_2^4(z) + \vartheta_4^4(z) \\
 \vartheta_1(y+z)\vartheta_1(y-z)\vartheta_4^2 & = & \vartheta_1^2(y) \vartheta_4^2(z) - \vartheta_4^2(y) \vartheta_1^2(z) \\
 \vartheta_4(2y) \vartheta_4^3 & = & \vartheta_4^4(y) - \vartheta_1^4(y) \\
 \vartheta_1(2y) \vartheta_2 \vartheta_3 \vartheta_4 & = & 2 \vartheta_1(y)\vartheta_2(y)\vartheta_3(y) \vartheta_4(y) \\
\end{eqnarray*}
The equation
$\vartheta_2^4 + \vartheta_4^4 = \vartheta_3^4$
has no non-trivial rational solution.
(Clearing denominators would lead to a solution of Fermat's equation for
exponent 4.)
This is in contrast to the situation with elliptic curves, where modular
solutions of the elliptic curve equation (``points on the curve'') can be turned into rational solutions of
an equivalent elliptic curve.

The restriction against rational solutions can be circumvented by
switching to another field where the equation has solutions.
Two obvious choices are an algebraic number field, or to work modulo
a prime number.
I decided to experiment with the modulus $P=43$.  A $4k+3$ prime was selected
so that more fourth powers exist and so that the number $i = \sqrt{-1}$ has independent meaning.
I selected trial values mod 43 for $\vartheta_{1-4}$ and $\vartheta_{1-4}(1)$.
The values were chosen to be compatible with the equations above.
Having selected these 8 numbers, we use theta function identities to compute a table of
$\vartheta_{1-4}$ evaluated at 2, 3, etc.

\begin{center}
\begin{tabular}{crrrr}
\multicolumn{5}{c}{Theta Functions Mod 43}\\
   &    &    &    &    \\
N & $\vartheta_1$ & $\vartheta_2$ & $\vartheta_3$ & $\vartheta_4$ \\
0  &  0 &  2 & 14 &  1 \\
1  &  1 & 11 &  7 &  2 \\
2  & 11 & 24 & 33 & 15 \\
3  &  1 &  1 & 24 & 24 \\
4  & 34 & 12 & 32 & 36 \\
5  & 17 & 25 &  7 &  3 \\
6  & 35 &  0 & 39 & 30 \\
7  & 38 & 20 &  3 & 32 \\
8  &  8 & 25 &  5 & 11 \\
9  & 16 & 27 & 40 & 40 \\
10 & 14 & 32 & 42 & 23 \\
11 & 24 & 37 & 39 &  5 \\
12 &  0 &  4 & 15 & 41 \\
13 &  7 & 34 & 37 & 29 \\
   &    &    &    &    \\
ratio 13/1 & 7 & 7 & -7 & -7 \\
   &    &    &    &    
\end{tabular}
\end{center}

Note that row 12 has a relationship to row 0, with $\vartheta_1$ and $\vartheta_2$ being doubled,
while $\vartheta_3$ and $\vartheta_4$ are multiplied by -2.
Rows 13 and 1 are similarly related, with ratios $\pm 7$.
We find a ratio relationship with consistent signs,
between rows 24 and 0, 25 and 1, etc.:
\vspace{1 em} \\
Row 24 = 16 * row 0;  row 25 = 24 * row 1;  row 26 = 36 * row 2.
\vspace{1 em} \\
With a large enough separation between rows,
  a constant ratio will develop,
and after some larger separation, a true period will appear.

An equation similar to the Somos recurrence works:
$$ \vartheta_1^2(n) = 8 \vartheta_1(n+1) \vartheta_1(n-1) + \vartheta_1(n+2) \vartheta_1(n-2).$$

We explore the ratio relationship between theta functions and elliptic curves:

$$ {\rm sn}(n) = \frac{\vartheta_3}{\vartheta_2} * \frac{\vartheta_1(n)}{\vartheta_4(n)} $$

\begin{tabular}{cll}
$n$   & sn($n$) & \\
0-11  & 0,25,8,20,9,11,1,11,9,20,8,25, &  this row is symmetric \\
12-23 & 0,18,35,23,...                 &  this row is the negative of the row above \\
 & &
\end{tabular} \\
Cn and dn work as well.
Checking out the elliptic curve that is parameterized by these modular elliptic functions,
$$ X= {\rm sn}(n); \;\;  Y = {\rm cn}(n) {\rm dn}(n) $$
$$ Y^2 = 1-7X^2+6X^4  = (1-X^2)(1-6X^2)  \;\;  {\rm ;elliptic \  curve \  equation} $$

Cn and sn are relabeled sine and cosine:
The set
$\{{\rm cn}(n) + i \, {\rm sn}(n)\}$  is the powers of $(13-2i) \bmod 43$, but in an apparently random order.
Note that $13-2i$ is on the mod 43 unit circle, since $||13-2i|| = 13^2 + 2^2 = 1 \bmod 43$.
The random ordering is captured by am, which in the real world is the distortion function:
$$
{\rm sn}(n) = {\rm sin}({\rm am}(n))
$$

\subsection{Some Possible Uses of Modular Theta Functions}

Modular theta functions may be directly useful in Diffie Hellman key exchange.
They could be used to compute elliptic curve values.
And they provide another example of a modular analog for a classical special function.

\section{Conclusion}

There's a lot more to learn here.

\bibliography{bib}

\end{document}